\documentclass[english]{elsarticle}
\usepackage[T1]{fontenc}
\usepackage[latin9]{inputenc}
\usepackage{color}
\usepackage{amsmath}
\usepackage{amssymb}
\usepackage{graphicx}

\makeatletter

\newcommand{\lyxdot}{.}

\usepackage{url}
\usepackage{fullpage}

\makeatother

\usepackage{babel}
\begin{document}
\title{An adaptive step size controller for iterative implicit methods\tnoteref{label1}}
\author[uibk]{Lukas Einkemmer\corref{cor1}} \ead{lukas.einkemmer@uibk.ac.at}
\address[uibk]{Department of Mathematics, University of Innsbruck, Austria}
\cortext[cor1]{Corresponding author}
\begin{abstract} 
The automatic selection of an appropriate time step size has been considered extensively in the literature. However, most of the strategies developed operate under the assumption that the computational cost (per time step) is independent of the step size. This assumption is reasonable for non-stiff ordinary differential equations and for partial differential equations where the linear systems of equations resulting from an implicit integrator are solved by direct methods. It is, however, usually not satisfied if iterative (for example, Krylov) methods are used.

In this paper, we propose a step size selection strategy that adaptively reduces the computational cost (per unit time step) as the simulation progresses, constraint by the tolerance specified. We show that the proposed approach yields significant improvements in performance for a range of problems (diffusion-advection equation, Burgers' equation with a reaction term, porous media equation, viscous Burgers' equation, Allen--Cahn equation, and the two-dimensional Brusselator system). While traditional step size controllers have emphasized a smooth sequence of time step sizes, we emphasize the exploration of different step sizes which necessitates relatively rapid changes in the step size.
\end{abstract}  
\maketitle

\section{Introduction}

Numerically solving time dependent differential equations is an important
task in many fields of science and engineering. Crucial to that process
is an efficient time integrator. Since the early advent of computers
such methods have been used to first solve ODEs (ordinary differential
equation) and then PDEs (partial differential equations).

Numerical simulations can be run with a constant time step size. However,
modern software packages usually automatically select an appropriate
step size given a desired tolerance (which is specified by the user).
To accomplish this so-called (automatic) step size controllers are
used in conjunction with an error estimator. Such an approach provides
a range of advantages. First, it frees the user from selecting an
appropriate step size and, ideally, from verifying the accuracy of
the simulation (by numerical convergence studies or similar means).
Second, a good step size controller is not only able to provide an
estimate of the error made, but also to detect the onset of numerical
instabilities and to reduce the step size to prevent them. This is
of particular importance for explicit methods which, for example,
can not operate with step sizes above the Courant\textendash Friedrichs\textendash Lewy
(CFL) limit and for implicit methods that are not A-stable. Last,
but certainly not least, step size controllers can increase the computational
efficiency by allowing the software to adaptively increase and decrease
the step size during the course of the simulation. This is usually
done in response to an error estimate, where errors significantly
below the specified tolerance indicate the possibility to increase
the time step.

Step size controllers require an error estimate. Fortunately, estimating
the error can often be accomplished with only a minor increase in
the computational cost. One approach commonly used are so-called embedded
Runge\textendash Kutta methods. These schemes consist of a pair of
Runge\textendash Kutta methods with different order that share most
or even all internal stages (see, for example, \cite[Chap. II.4]{hairerI}
or \cite{cash1979,dormand1980,shampine1984}). For multistep methods
a comparison with extrapolated values is often used (see, for example,
\cite{cvode}). Alternatively, Richardson extrapolation can be used
to obtain an error estimate, but is usually more demanding from a
computational point of view.

Almost all step size controllers are predicated on the assumption
that the largest possible step size should be selected. Thus, the
step size is chosen such that the error committed exactly matches
the tolerance specified by the user. This is a reasonable assumption
for explicit Runge\textendash Kutta methods, where the computational
cost is independent of the step size. Now, let us assume that our
error estimator provides an estimate $\epsilon^{k}$ for the $k$th
step (note that in accordance with much of the PDE literature we use
superscripts to denote the time indices). The local error of a numerical
method with order $p$ is modeled as $e^{k}=D(\tau^{k})^{p+1}$, where
$\tau$ denotes the time step size and $D$ is a constant (which for
the purpose of simplicity is assumed to be independent of $k$; in
most situations this is a reasonable assumption as the error constant
only varies slowly on $\mathcal{O}(\tau)$ timescales). Then, to determine
the optimal step size we set
\[
{\normalcolor {\color{blue}{\color{black}\text{tol}=e^{k+1}=D(\tau^{k+1})^{p+1},}}}
\]
where $\text{tol}$ is the user specified tolerance. This is not a
particular useful constraint to determine $\tau^{k+1}$ as $D$ is
unknown. Thus, we consider
\[
\frac{e^{k+1}}{e^{k}}=\left(\frac{\tau^{k+1}}{\tau^{k}}\right)^{p+1}
\]
which can be solved for $\tau^{k+1}$
\begin{equation}
\tau^{k+1}=\tau^{k}\left(\frac{\text{tol}}{e^{k}}\right)^{1/p}.\label{eq:P-controller}
\end{equation}
Equation (\ref{eq:P-controller}) allows us to estimate the optimal
time step $\tau^{k+1}$ based on the previous error estimate $e^{k}$
a\textcolor{black}{nd the previous time step size $\tau^{k}$. This
then results in a sequence of times $t^{k}$ at which a numerical
approximation is obtained. In pra}ctice this is a dangerous approach
as even very small errors in the error estimate can result in a time
step sizes that exceeds the prescribed tolerance (leading to step
size rejection). Thus, usually a safety factor is incorporated. For
more details we refer the reader to \cite[Chap. II.4]{hairerI} and
\cite{shampine2005}.

This simple formula can be interpreted as a P controller. The mathematical
analysis is in fact based on this observation (see, for example, \cite{gustafsson1988,gustafsson1994,soderlind2002,soderlind2006}).
Consequently, PI controllers have been introduced \cite{gustafsson1988},
which for some integrators and problems show an increase in performance.
Certainly, these PI controllers increase the smoothness of the step
size sequence (i.e.~the change in step size behaves less erratic).
These ideas have been enhanced in a variety of directions. The importance
of changing strategies when operating close to the stability limit
for explicit methods has also been recognized \cite{hall1995}.

Although some work has been conducted in estimating global errors
(see, for example, \cite{shampine2005}), the local step size controllers
described above, with some modifications to avoid excessively large
step sizes, still form the backbone of most time integration packages.
For example, the RADAU5 code \cite[Chap. IV.8]{hairerII} employs
a variant of the PI controller, while the multistep based CVODE code
\cite{cvode} and \cite{eckert2004} uses a variant of the P controller.
As a result, the described step size controllers have been extensively
tested and used in a range of applications, both for ordinary as well
as for partial differential equations.

The desire for solving partial differential equations with ever increasing
grid sizes and more accurate physical models, however, calls into
question the validity of the assumptions made. In both the RADAU5
and CVODE code mentioned above implicit numerical methods are employed
to solve the stiff ODEs resulting from the space discretization of
the PDE under consideration. These implicit methods require the solution
of a linear system of equations which is now routinely done by iterative
numerical methods (such as the conjugate gradient method or GMRES).
However, the number of iterations required is quite sensitive to the
linear system solved. In particular, smaller time step sizes reduce
the magnitude of the largest eigenvalue of the matrix, which in turn
reduces the number of iterations required per time step. This means
that reducing the time step size below what is dictated by the specified
tolerance, e.g.~according to equation (\ref{eq:P-controller}), can
actually result in an increase in performance.

\textcolor{black}{Many implementations do not exploit this fact. However,
the issue at hand has been recognized in \cite{hochbruck1998} and
\cite{weiner1998}. Both of these approaches limit the size of the
Krylov subspace. In \cite{hochbruck1998} both lower and upper bounds
are specified. If the Krylov dimension falls within those bounds,
the time step size is chosen according to the step size controller.
If this is not the case, the step size is adjusted. In \cite{weiner1998}
a multiple Arnoldi process is used. In that context the increase of
the Krylov dimension in the higher stages is limited by a fixed value.
The downside of this approach is that it is usually not known a priori
how the bounds should be chosen; the corresponding value is most likely
highly problem dependent. It sho}uld also be noted that in the context
of ODEs the importance of considering variations in cost (as a function
of the time step size) has been recognized in \cite{gustafsson1997}.
There are analytically derived estimates of the cost and, in line
with the control theoretic approach to step size selection, smooth
step size sequences have been emphasized. In contrast, in this work
we emphasize dynamically obtained estimates of the cost (which is
particularly useful for nonlinear PDEs, where obtaining good a priori
estimates is often extremely difficult) and exploration of the space
of admissible step sizes (which results in rather frequent and often
erratic step size changes). Furthermore, in the context of approximating
matrix exponentials by polynomial interpolation at Leja points, a
procedure to determine the optimal step size based on a backward error
analysis has been proposed \cite{caliari2016}. This approach can
be very effective but requires certain information on the spectrum
of the matrix under consideration. This information is not easily
obtained in a matrix free implementation and, for nonlinear PDEs,
can change from one time step to the next. In contrast, the step size
controller proposed in this work requires no a priori information
and is thus designed to naturally work for matrix free implementations.

In addition to the considerations above, it has been observed in many
applications that a reversed C shape can be observed on a precision-work
diagram for the traditional step size controllers. That is, specifying
a more stringent tolerance initially results in an increase of performance
(i.e.~smaller run times). The problem with that approach is that
the user of the software is once again tasked with finding the best
step size (or rather decreasing the tolerance until the run time is
minimized as well). Thus, effectively counteracting one of the primary
advantages of automatic step size control.

Such behavior can be observed across a range of test problems \cite{loffeld2013,luan2017,hochbruck1998}
as well as for problems that stem from more realistic physical models
\cite{einkemmer2017,blom2016,narayanamurthi2017}. As the before mentioned
work shows, this behavior is not limited to one class of numerical
method but can be observed for implicit Runge\textendash Kutta methods,
BDF methods, implicit-explicit (IMEX) methods, and exponential integrators. 

In this paper we propose an approach for adaptive step size control
that does not optimize for the largest time step size but rather tries
to minimize computational cost. Since it is difficult to analytically
determine the optimal step size, an optimization procedure is used
in parallel with the time stepper. The basic idea of this algorithm
is described in section \ref{sec:Basic-algorithm}. There numerical
results for a linear diffusion-advection equation are also shown.
In section \ref{sec:Nonlinear-problems} the efficiency of our algorithm
applied to four nonlinear problems is investigated. Finally, we conclude
in section \ref{sec:Conclusion}.

\section{Basic algorithm\label{sec:Basic-algorithm}}

\subsection{Setting\label{subsec:Setting}}

For the remainder of this section we consider the linear diffusion-advection
equation
\begin{equation}
\partial_{t}u(t,x)=\partial_{xx}u(t,x)+\eta\partial_{x}u(t,x)\label{eq:diffadv}
\end{equation}
with periodic boundary conditions on $[0,1]$. The dimensionless Péclet
number $\eta$ determines the relative strength of advection compared
to diffusion. As initial value we prescribe the following Gaussian
\[
u(0,x)=\mathrm{e}^{-(x-1/2)^{2}/(2\sigma_{0}^{2})}.
\]
\textcolor{black}{In this setting the analytic solution is known exactly
(strictly speaking this is only true if the problem is posed on the
entire real line; however, for times where the spread of the Gaussian
is smaller than the computational domain, a similar dynamics can be
observed for periodic boundary conditions). The Gaussian is translated
in space and the bump spreads out (assuming that the solution is sufficiently
small at the boundary). The spread of the standard deviation $\sigma$
at time $t$ is given by $\sigma(t)=\sqrt{\sigma_{0}^{2}+2t}$. In
particular, this implies that if a small $\sigma_{0}$ is chosen the
time step size is initially dictated by accuracy constraints (even
for small Péclet numbers). However, later in the evolution implicit
time integrators can take relatively large time steps without incurring
a significant error. In the numerical simulations conducted we will
choose a final time $t=0.2$. Thus, the present test problem probes
both of these regimes.}

\textcolor{black}{Note that in the numerical simulations we will present,
advection dominated (i.e.~large $\eta$) processes are considered
as well. In this regime, an explicit numerical methods could also
be used. However, it is our view that a general purpose implicit integrator/step
size controller should also be able to handle such problems. We will
see, however, that this is a challenging problem for the traditional
step size controller.}

\textcolor{black}{In all our implementations we use the standard centered
difference schem}e to discretize the diffusive part and a simple upwind
scheme for the advection. Concerning the time discretization, implicit
Runge\textendash Kutta methods will be employed. This is done in order
to avoid some of the tedious details encountered when dealing with
variable step size multistep methods (for example, limitations on
how rapidly the time step size is allowed to change). 

Many implicit Runge\textendash Kutta method have been considered in
the literature. Perhaps, the most well known are the Crank\textendash Nicolson
method and the classes of Gauss and Radau methods. The latter forms
the basis for the widely used RADAU5 time integrator. However, the
issue with higher order collocation methods is that a straightforward
implementation yields large matrices to invert. This further worsens
the problem that we try to overcome in this paper (see below) which
incidentally would lead to overly optimistic results. In addition,
almost all of the publicly available integrator packages either do
not directly support sparse Krylov solvers (such as RADAU5), make
it very difficult to change the time stepping strategy or to specify
a fixed order or fixed time step size (such as CVODE). Thus, in the
following we will use the Crank\textendash Nicolson scheme along with
a two stage third order SDIRK (singly diagonally implicit) scheme,
henceforth called SDIRK23, and a five stage fourth order SDIRK scheme,
henceforth called SDIRK54. All of these methods require us to (only)
solve an $n\times n$ linear system, where $n$ is the number of grid
points and thus give a better indication of the actual performance
attainable. The same can be accomplished for an implementation of
Radau methods (see, for example, \cite[Chap. IV.8]{hairerII}) but
in this case the details of the implementation are much more involved.
This is a further reason to stick with the, relatively simple, Crank\textendash Nicolson,
SDIRK23, and SDIRK54 method.

Now, let us describe the numerical time integrators used in this paper
in more detail. After an appropriate discretization in space we have
to integrate the following system of ODEs (ordinary differential equations)
\[
y^{\prime}(t)=f(t,y(t))
\]
in time. For linear autonomous problems this ODE could be significantly
simplified. However, since we will encounter nonlinear problems in
section \ref{sec:Nonlinear-problems} we will consider the more general
formulation here. In this setting the Crank\textendash Nicolson scheme
is given by
\[
y^{1}=y^{0}+\frac{\tau}{2}\left(f(0,y^{0})+f(\tau,y^{1})\right),
\]
where the time step size of $\tau$ is conducted to obtain $y^{1}$
from $y^{0}$. The third order SDIRK23 scheme is given by (see, for
example, \cite{ascher1997})

\begin{align*}
k^{1} & =f(\gamma\tau,y^{0}+\tau\gamma k^{1})\\
k^{2} & =f((1-\gamma)\tau,y^{0}+(1-2\gamma)\tau k^{1}+\gamma\tau k^{2})\\
y^{1} & =y^{0}+\frac{\tau}{2}\left(k^{1}+k^{2}\right),
\end{align*}
where $\gamma=\frac{3+\sqrt{3}}{6}$. For a linear problem (i.e.~$f(t,y(t))=Ay(t)$
for $A\in\mathbb{R}^{n\times n}$) this yields
\begin{align*}
(I-\tau\gamma A)k^{1} & =Ay^{0}\\
(I-\tau\gamma A)k^{2} & =Ay^{0}+(1-2\gamma)\tau Ak^{1}\\
y^{1} & =y^{0}+\frac{\tau}{2}(k^{1}+k^{2}).
\end{align*}
Thus, we only have to solve two $n\times n$ linear systems. In the
nonlinear case this is still true but the linear solve is now conducted
as the inner loop in Newton's method. The SDIRK54 scheme is given
by \cite[p. 107]{hairerII}
\begin{align*}
k^{1} & =\tau f\left(\tfrac{1}{4}\tau,y^{0}+\tfrac{1}{4}k^{1}\right)\\
k^{2} & =\tau f\left(\tfrac{3}{4}\tau,y^{0}+\tfrac{1}{2}k^{1}+\tfrac{1}{4}k^{2}\right)\\
k^{3} & =\tau f\left(\tfrac{11}{20}\tau,y^{0}+\tfrac{17}{50}k^{1}-\tfrac{1}{25}k^{2}+\tfrac{1}{4}k^{3}\right)\\
k^{3} & =\tau f\left(\tfrac{1}{2}\tau,y^{0}+\tfrac{371}{1360}k^{1}-\tfrac{137}{2720}k^{2}+\tfrac{15}{544}k^{3}+\tfrac{1}{4}k^{4}\right)\\
k^{5} & =\tau f\left(\tau,y^{0}+\tfrac{25}{24}k^{1}-\tfrac{49}{48}k^{2}+\tfrac{125}{16}k^{3}-\tfrac{85}{12}k^{4}+\tfrac{1}{4}k^{5}\right)\\
y^{1} & =y^{0}+\tfrac{25}{24}k^{1}-\tfrac{49}{48}k^{2}+\tfrac{125}{16}k^{3}-\tfrac{85}{12}k^{4}+\tfrac{1}{4}k^{5}.
\end{align*}
This SDIRK54 method is L-stable. The Crank\textendash Nicolson method
and the SDIRK23 method are A-stable but not L-stable. Even though
we use A-stable implicit methods here it is instructive to consider
the stability constraints that the explicit Euler method would encounter.
In this case we have
\[
\tau<\min\left(\frac{1}{2}\frac{1}{n^{2}},\frac{1}{\vert\eta\vert n}\right),
\]
where $n$, as before, denote the number of grid points. Usually,
this is dominated by the stability constraint from the diffusion (the
first term in the formula) but we will also consider examples where
the Péclet number is large enough such that these two stability constraints
are comparable.

In all our examples we will use the GMRES (generalized minimal residual)
method to solve the resulting linear system (note that for $\eta\neq0$
the matrix $A$ is not symmetric). This iterative Krylov subspace
method terminates when the residual is below one-tenth of the tolerance
prescribed for the numerical method. It should be emphasized, however,
that the step size controller proposed in the next section is completely
independent of the iterative method used. The choice of Krylov methods
is due to their ubiquity in applications. Nonetheless, relaxation
methods or methods based on (direct) polynomial interpolation could
be used just as easily (as all the relevant data are obtained at run
time).

\subsection{Step size controller\label{subsec:Step-size-controller}}

As has been outlined in the introduction, traditional step size controller
always take the largest step possible given the accuracy constraints.
In our case, we will adaptively change the step size depending on
the cost of the previous time step. This allows us to explore a range
of step sizes and search for the most cost effective one (which might
be significantly smaller than the one selected by a traditional step
size controller). In the following we will interpret this as a one-dimensional
gradient descent optimization algorithm.

\textcolor{black}{Our goal is to optimize the computational cost per
unit time step, i.e.
\[
{\color{blue}{\color{black}c^{k}=\frac{M(i^{k})}{\tau^{k}},}}
\]
where $i^{k}$ is the number of Krylov iterations conducted in the
$k$th time step (i.e.~the sum of the Krylov dimensions over all
stages of the method) and $\tau^{k}$ is the size of the step. The
function $M\colon\mathbb{N}\to\mathbb{R}_{\geq0}$ models the computational
cost as a function of the number of iterations. This is not a trivial
task and the function $M$, in general, depend on the computer system
used. For example, on large scale supercomputers the latency introduced
by the dot products can actually be the limiting factor with respect
to performance. The situation is further complicated by the fact that
GMRES has to be used with a restart procedure (in our simulation we
will restart every $20$th iterations). In the present work we will
assume that the computational cost is directly proportional to the
number of iterations; i.e.~$M(i^{k})=i^{k}$. This assumption is
valid if the cost of the dot products is small compared to the cost
of evaluating the right-hand side of the PDE.  However, we duly note
that the controller presented in the following can be used just as
well with a different cost model $M$ and that there are certainly
situations where this would be indicated.}

\textcolor{black}{Since the number of Krylov iterations is determined
adaptively, th}e cost is only available at the end of the time step.
Our goal is to dynamically adjust the step size (i.e.~to change $\tau^{k}$)
such that 
\[
c^{k}\to\text{min}.
\]
As is usually done in the analysis of step size controllers, we work
with the logarithm of the step size $T^{k}=\log\tau^{k}$ and the
logarithm of the computational cost $C^{k}(T^{k})=\log c^{k}(\tau^{k})$.
Employing a one-dimensional gradient descent algorithm we have (assuming
that the cost per unit step $C^{k}$ is only a function of the time
step size $T^{k}$)
\[
T^{k+1}=T^{k}-\gamma\nabla C^{k}(T^{k})
\]
Now, since we have only discrete values at our disposal we approximate
the gradient by a difference quotient
\[
\nabla C^{k}(T^{k})\approx\frac{C^{k}(T^{k})-C^{k}(T^{k-1})}{T^{k}-T^{k-1}}.
\]
By doing so we effectively rule out taking a constant time step size
(i.e.~$T^{k}=T^{k-1}$). This is necessary in order to obtain the
necessary information that guide our adaptive step size selection.
In the literature it is often argued that a smooth step size selection
is desirable in order to increase the accuracy of the error estimator.
However, in our case we will only use the error estimator as a worst
case bound. Most of the time, the step size will be chosen well below
that limit. This implies that step size rejection happens infrequently
(if at all) even if the time step size varies considerably from one
step to the next.

Unfortunately, making this approximation is not yet sufficient as,
strictly speaking, we have to make a distinction between the cost
functions $C^{k-1}$ and $C^{k}$. Note that $C^{k-1}(T)$ gives the
cost of making a step with size $T$ starting from the beginning of
the previous time step, i.e.~starting from $t^{k-1}$. This value
is, in general, different from $C^{k}(T)$ (the cost of making a step
with the same size but starting at $t^{k}$). During the time integration
only $C^{k-1}(T^{k-1})$ but not $C^{k}(T^{k-1})$ is sampled. Thus,
we write
\[
\frac{C^{k}(T^{k})-C^{k}(T^{k-1})}{T^{k}-T^{k-1}}=\frac{C^{k}(T^{k})-C^{k-1}(T^{k-1})}{T^{k}-T^{k-1}}+\frac{C^{k-1}(T^{k-1})-C^{k}(T^{k-1})}{T^{k}-T^{k-1}}.
\]
Further assuming that $C^{k}$ changes slowly as a function of $k$
we obtain
\[
\nabla C^{k}(T^{k})\approx\frac{C^{k}(T^{k})-C^{k-1}(T^{k-1})}{T^{k}-T^{k-1}}.
\]
This then gives us
\[
T^{k+1}=T^{k}-\gamma\frac{C^{k}(T^{k})-C^{k-1}(T^{k-1})}{T^{k}-T^{k-1}}
\]
In principle $\gamma$ is a free parameter (which can depend on both
$c$ and $\tau$). Taking the exponential on both sides gives
\begin{equation}
\tau^{k+1}=\tau^{k}\exp\left(-\gamma\Delta\right)\label{eq:tauk+1generic}
\end{equation}
with 
\[
\Delta=\frac{\log c^{k}-\log c^{k-1}}{\log\tau^{k}-\log\tau^{k-1}}.
\]
The simplest choice is taking $\gamma=\text{const}$ which, however,
has two major drawbacks. First, if the cost varies rapidly for a relatively
small change in $\tau$ we obtain extremely large time step changes.
Second, if the cost varies slowly we might only change the time step
very slowly which means we can not explore the available parameter
space efficiently. In fact, as we will see later, this ability to
change the time step size even if $\Delta$ is small will be a crucial
ingredient of our step size controller. Therefore, we propose to use
the following method to determine the new time step size $\overline{\tau}^{k+1}$
from the old time step size $\tau^{k}$ 
\begin{align}
s & =\exp\left(-\alpha\tanh\left(\beta\Delta\right)\right)\nonumber \\
\overline{\tau}^{k+1} & =\tau^{k}\begin{cases}
\lambda & 1\leq s<\lambda\\
\delta & \delta\leq s<1\\
s & \text{otherwise}
\end{cases}\label{eq:tauk+1}
\end{align}
which can be realized by choosing $\gamma(\Delta)$ in equation (\ref{eq:tauk+1generic}).
This ensures that the time step size is changed by at least $\lambda\tau^{k}$
or $\delta\tau^{k}$ and limits the maximum step size change to $\exp(\pm\alpha)\tau^{k}$.
The understanding here is that the parameters $\alpha,\beta,\delta$
and $\lambda$ are positive and that $\lambda>1$ and $\delta<1$
are sufficiently separated from $1$ in order to always cause a non-trivial
change in step size. The new time step size $\overline{\tau}^{k+1}$
is used only if it is smaller or equal to the time step size determined
from the traditional controller. This is necessary as there is no
guarantee that the step size computed from equation (\ref{eq:tauk+1})
satisfies the accuracy requirement specified by the user. That is,
in all simulations we employ
\begin{equation}
\tau^{k+1}=\text{min}\left(\overline{\tau}^{k+1},\tau^{k}\left(\frac{\text{tol}}{e^{k}}\right)^{1/p}\right).\label{eq:proposed-controller}
\end{equation}

Limiting ourselves to a functional representation with four free parameters
implies that we obtain a manageable optimizing problem for these parameters.
The goal function is set by considering the average performance for
$(n,\eta)\in\mathcal{N}=\{(100,10),(300,100),(500,0),(500,1000)\}$
and $\text{\ensuremath{\epsilon}}\in\mathcal{E}=\{10^{-2},10^{-3},10^{-4},10^{-5},10^{-7}\}$,
where $\epsilon$ is the tolerance specified for the numerical method.
We integrate until final time $t=0.2$ and use the SDIRK54 method
with $\sigma_{0}=1.4\cdot10^{-3}$ for the initial value. We consider
the two fitness functions $f_{1}(\Gamma)$ and $f_{2}(\Gamma)$, depending
on the parameter $\Gamma=(\alpha,\beta,\lambda,\delta)$, given by
\[
f_{1}(\Gamma)=\sum_{(n,\eta)\in\mathcal{N}}\sum_{\epsilon\in\mathcal{E}}R_{n\eta\epsilon}(\Gamma)
\]
and 
\[
f_{2}(\Gamma)=\sum_{(n,\eta)\in\mathcal{N}}\sum_{\epsilon\in\mathcal{E}}Q_{\Pi}(R_{n\eta\epsilon}(\Gamma),T_{n\eta\epsilon}),
\]
where $R_{n\eta\epsilon}(\Gamma)$ is the run time of our step size
controller with parameter $\Gamma$ for the linear diffusion-advection
problem specified by $(n,\eta,\epsilon)$, $T_{n\eta\epsilon}$ is
the run time required by the traditional step size controller for
the same problem, and
\[
Q_{\Pi}(x,y)=\frac{x}{y}\cdot\begin{cases}
1 & x/y\leq1\\
\Pi & \text{otherwise}
\end{cases}.
\]
The constant parameter $\Pi$ penalizes cases where our approach performs
worse than the traditional step size controller. This penalty parameter
is used to trade off increased gain in performance (in, hopefully,
the majority of configurations) with how much reduction in performance
we are willing to tolerate (for, hopefully, a small number of configurations).

\begin{figure}
\centering{}\includegraphics[width=14cm]{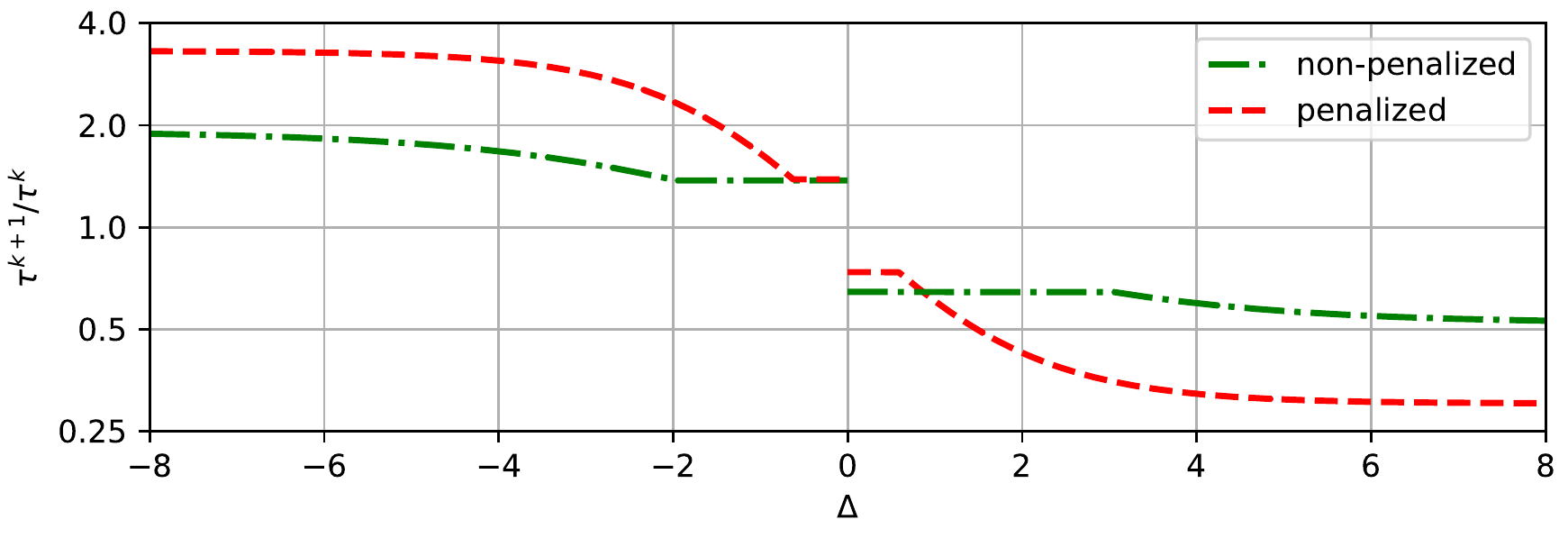}\caption{The ratio of the new step size to the old step size (i.e.~$\tau^{k+1}/\tau^{k}$)
is shown as a function of $\Delta$ for both the non-penalized (green
dash-dotted line) and the penalized step size controller (red dashed
line). \label{fig:optimized}}
\end{figure}

In the following we will consider two step size controllers based
on $f_{1}(\Gamma)$ and $f_{2}(\Gamma)$ with penalty parameter $\Pi=10$.
A numerical optimization using differential evolution is performed.
For the former case this results in the set of parameters
\[
\alpha=0.65241444,\qquad\beta=0.26862269,\qquad\lambda=1.37412002,\qquad\delta=0.64446017.
\]
while the latter case gives

\[
\alpha=1.19735982,\qquad\beta=0.44611854,\qquad\lambda=1.38440318,\qquad\delta=0.73715227,
\]
We will refer to these two sets of parameters (and to the corresponding
step size controller) as the non-penalized and the penalized controller,
respectively. The corresponding curves given by (\ref{eq:tauk+1})
are shown in Figure \ref{fig:optimized}. The difference between these
two can be understood as follows. For the penalized integrator it
is often favorable to relatively rapidly increase the time step size
in order to quickly reach the upper limit dictated by the accuracy
requirement. If this regime is reached, the penalized step size controller
effectively acts like the traditional controller (with the difference
that it can still decreases the step size repeatedly in order to check
if a smaller step size reduces the computational cost). For the non-penalized
controller the goal is to more closely follow the gradient to find
the (locally) most efficient time step size. This means that it has
less variation in the time step size. 

It can also be observed from Figure \ref{fig:optimized} that the
minimal change in the step size is at least 30\%. This results in
the discontinuity seen in the plot at $\Delta=0$ and supports our
assertion that the step size should be varied rapidly in order to
collect the data necessary for our algorithm.

In this work we will compare the step size controller proposed in
this section with the traditional P controller. We choose the P controller
(as opposed to a PI controller or any of the other extensions that
have been developed) for two primary reasons. First, the P controller
is still widely used in software packages such as CVODE \cite{cvode}
and the advantage more elaborate schemes can provide depends on the
numerical method and a range of other factors. For example, in \cite{usman2000}
it shown that the PI controller fails to give good results for certain
multistep methods. Second, the P controller is a simpler algorithm
that only includes a single parameter (the safety factor) and most
of the results achieved carry over easily to more elaborate methods.

\subsection{Numerical results}

\begin{figure}
\begin{centering}
\includegraphics[width=16cm]{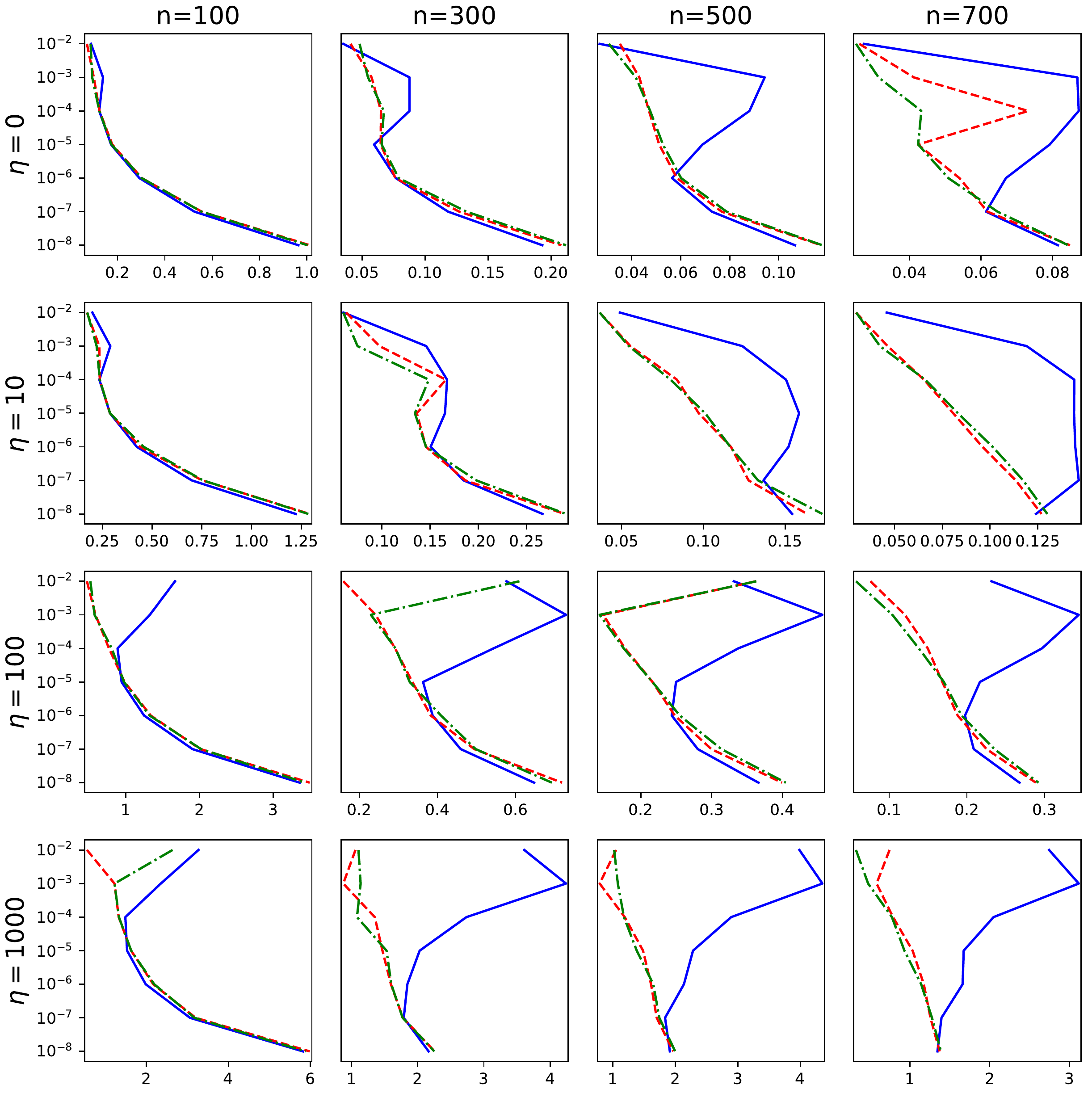}
\par\end{centering}
\caption{The number of normalized Krylov iterations employed by the Crank\textendash Nicolson
scheme (i.e.~computational cost; on the $x$-axis) for a given tolerance
(on the $y$-axis) is shown for equation (\ref{eq:diffadv}) (linear
diffusion-advection equation). The proposed step size controllers
are shown as dashed red lines (penalized variant) and dash-dotted
green lines (non-penalized variant), while the traditional step size
controller is shown in solid blue. The grid size ($n$ is the number
of grid points) and the Péclet number $\eta$ are varied.\label{fig:diffadv-cn}}
\end{figure}

We will use exactly the setup described in section \ref{subsec:Setting}.
For the initial value we choose $\sigma_{0}=1.4\cdot10^{-3}$ and
integrate until time $t=0.2$. In Figure \ref{fig:diffadv-cn} the
work-precision diagram for the Crank\textendash Nicolson method is
shown for the linear diffusion-advection equation. The cost (on the
$x$-axis) is represented as the number of Krylov iterations required
to advance the numerical solution by the same amount as the maximal
step size allowed by the classic CFL condition (i.e.~the number of
Krylov iterations is normalized to the cost of an explicit Euler method
that is operated with unit CFL number). The tolerance specified by
the user is shown on the $y$-axis. The blue line for the traditional
step size controller (a P controller with safety factor $0.9$) shows
the characteristics reversed C shaped curve mentioned in the introduction.
Ideally, we would expect from a step size controller a monotonous
increase in the cost as the tolerance is increased. We can see from
\ref{fig:diffadv-cn} that the step size controller proposed in section
\ref{subsec:Step-size-controller} (the dashed red and dash-dotted
blue line in the figure) matches this pattern very well (for a range
of grids and different Péclet numbers). Instead of the reversed C
shaped curve we now have, in almost all cases, a monotonously increase
of cost as a function of tolerance. What is perhaps even more important
is the significant decrease in computational cost that is obtained
by straightening out the reversed C curve. The actual increase in
computational performance varies with the specific configuration but
can yield a speedup of up to a factor of four. We also see that the
non-penalized variants gives generally better results.

\begin{figure}
\centering{}\includegraphics[width=16cm]{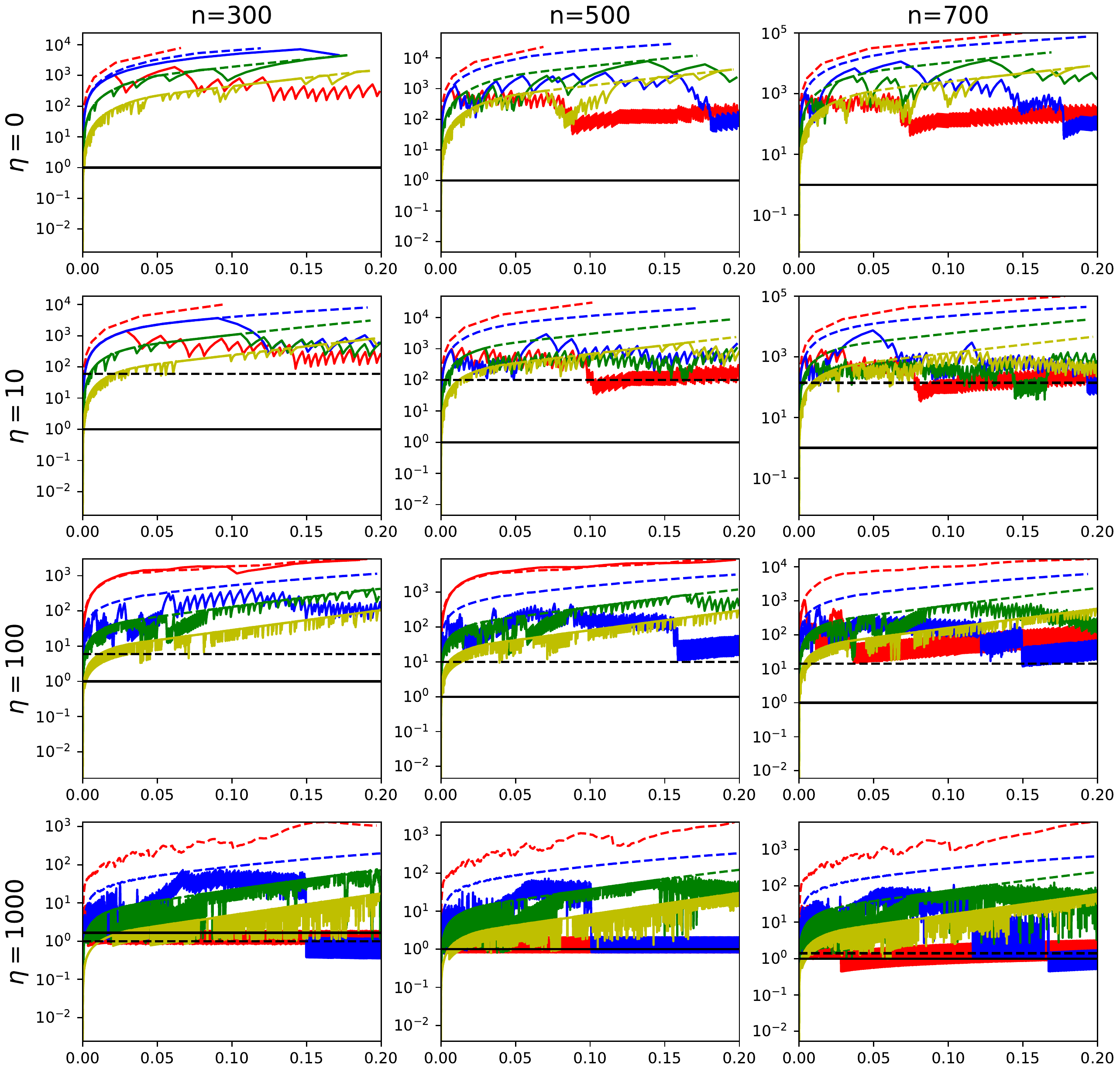}\caption{The CFL number (i.e.~size of the time step divided by the maximal
stable time step for the Euler method) taken by the Crank\textendash Nicolson
scheme is shown as a function of time (for the linear diffusion-advection
equation (\ref{eq:diffadv})). The solid lines correspond to the proposed
step size controller (the non-penalized variant), while dashed lines
correspond to the classic step size controller. In both cases results
for the tolerance set to $10^{-2}$ (red), $10^{-4}$ (blue), $10^{-6}$
(green), and $10^{-8}$ (yellow) are shown. The black line shows the
CFL condition induced by the diffusion and the dashed black line shows
the CFL condition induced by the advection.\textbf{ }The grid size
($n$ is the number of grid points) and the Péclet number $\eta$
are varied.\label{fig:diffadv-cn-hlist}}
\end{figure}

\begin{figure}
\begin{centering}
\includegraphics[width=16cm]{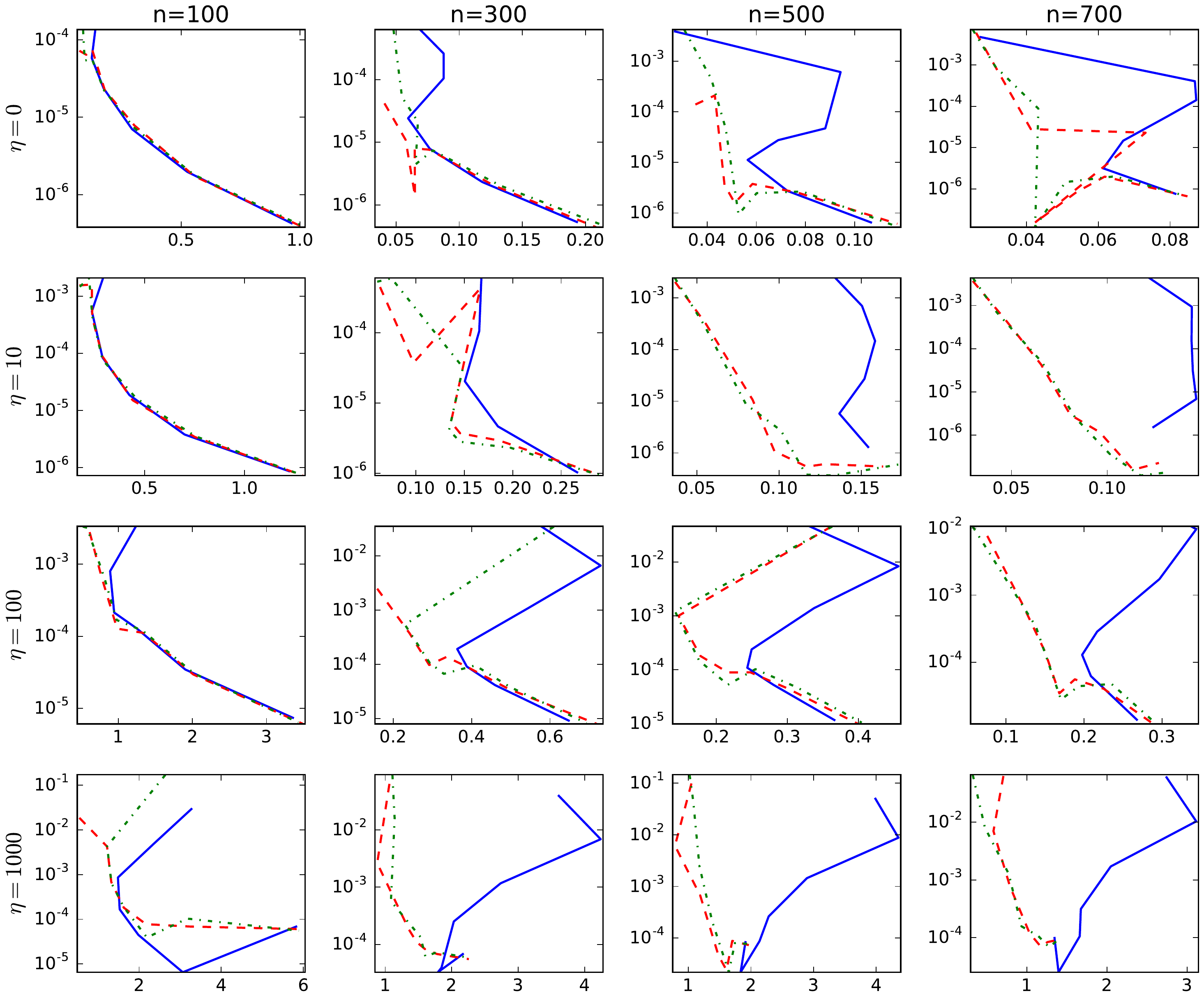}
\par\end{centering}
\caption{The number of normalized Krylov iterations employed by the Crank\textendash Nicolson
scheme (i.e.~computational cost; on the $x$-axis) for the achieved
accuracy at final time $t=0.2$ (on the $y$-axis) is shown for equation
(\ref{eq:diffadv}) (linear diffusion-advection equation). The proposed
step size controllers are shown as dashed red lines (penalized variant)
and dash-dotted green lines (non-penalized variant), while the traditional
step size controller is shown in solid blue. The grid size ($n$ is
the number of grid points) and the Péclet number $\eta$ are varied.\label{fig:diffadv-cn-acc}}
\end{figure}

Let us discuss this results in more detail. In Figure \ref{fig:diffadv-cn-hlist}
the step size for both the traditional controller (dashed lines) and
the proposed controller (the non-penalized variant) are shown for
four different tolerances. As has been discussed in the previous section,
the maximal allowable time step (due to accuracy constraints) increases
with time. This is in fact, what we observe for the traditional step
size controller. In the beginning of the simulation this curve is
very closely followed by our step size controller as well. This is
to be expected as in this case the integrator operates very close
to the accuracy limit. Later on, however, our step size controller
only increases the time step size if this results in reduced computational
cost (independent of accuracy constraints). This results in significantly
smaller time steps which in turn results in a significant reduction
in computational cost. This has the additional benefit that the numerical
solution is more accurate than the tolerance requested indicates (see
Figure \ref{fig:diffadv-cn-acc} and the following discussion). Figure
\ref{fig:diffadv-cn-hlist} also shows the that our algorithm dictates
a rapid change of step sizes. This is in stark contrast to the work
on control theoretic step size selection, where it was often argued
that a good control system should provide a smooth response \cite{gustafsson1988,soderlind2006}.
In our scheme, however, this rapid change is an important feature
in order to collect the necessary information to guide step size selection
(as has been discussed in the previous section). 

The results in Figure \ref{fig:diffadv-cn} do not provide information
on the actual global error that is achieved in the simulation. The
behavior of the global error can, of course, be significantly different
from the local error. The behavior of the corresponding error propagation,
for example, depends on the problem under consideration. However,
since the proposed step size controller is only able to decreases
the time step size compared to the traditional approach, the local
error per unit time step committed by the time integration scheme
is also reduced. This then implies that the global error of the proposed
controller is similar or smaller than the global error for the traditional
approach. To illustrate this we plot, in Figure \ref{fig:diffadv-cn-acc},
the global error at the final time as a function of the Krylov iterations.
If we compare these results to Figure \ref{fig:diffadv-cn} we can
see that the advantage of the proposed step size controllers is even
more pronounced according to that metric.

\begin{figure}
\begin{centering}
\includegraphics[width=16cm]{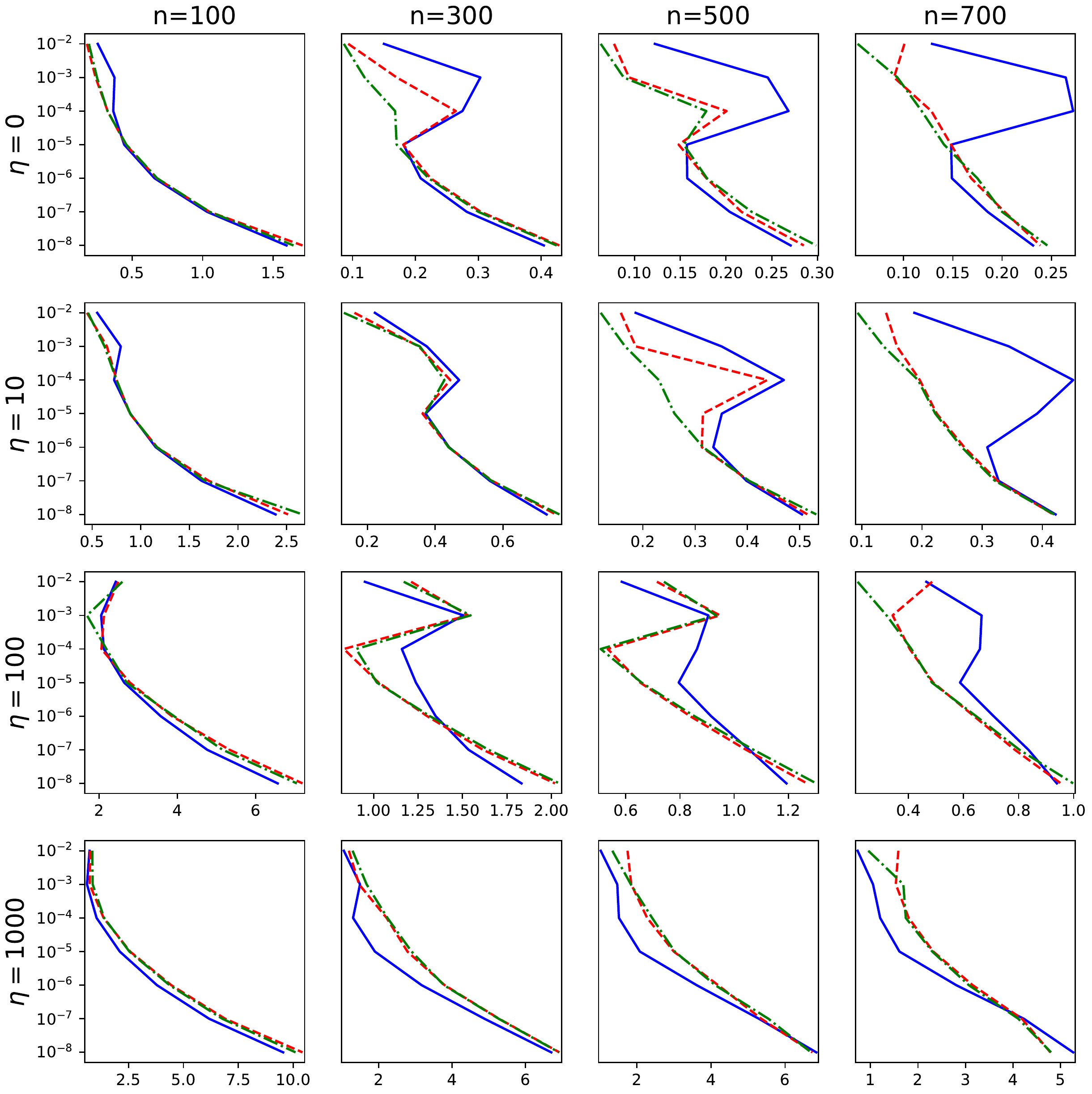}
\par\end{centering}
\caption{The number of normalized Krylov iterations employed by the SDIRK23
scheme (i.e.~computational cost; on the $x$-axis) for a given tolerance
(on the $y$-axis) is shown for equation (\ref{eq:diffadv}) (linear
diffusion-advection equation). The proposed step size controllers
are shown as dashed red lines (penalized variant) and dash-dotted
green lines (non-penalized variant), while the traditional step size
controller is shown in solid blue. The grid size ($n$ is the number
of grid points) and the Péclet number $\eta$ are varied.\label{fig:diffadv-sdirk23}}
\end{figure}

\begin{figure}
\begin{centering}
\includegraphics[width=16cm]{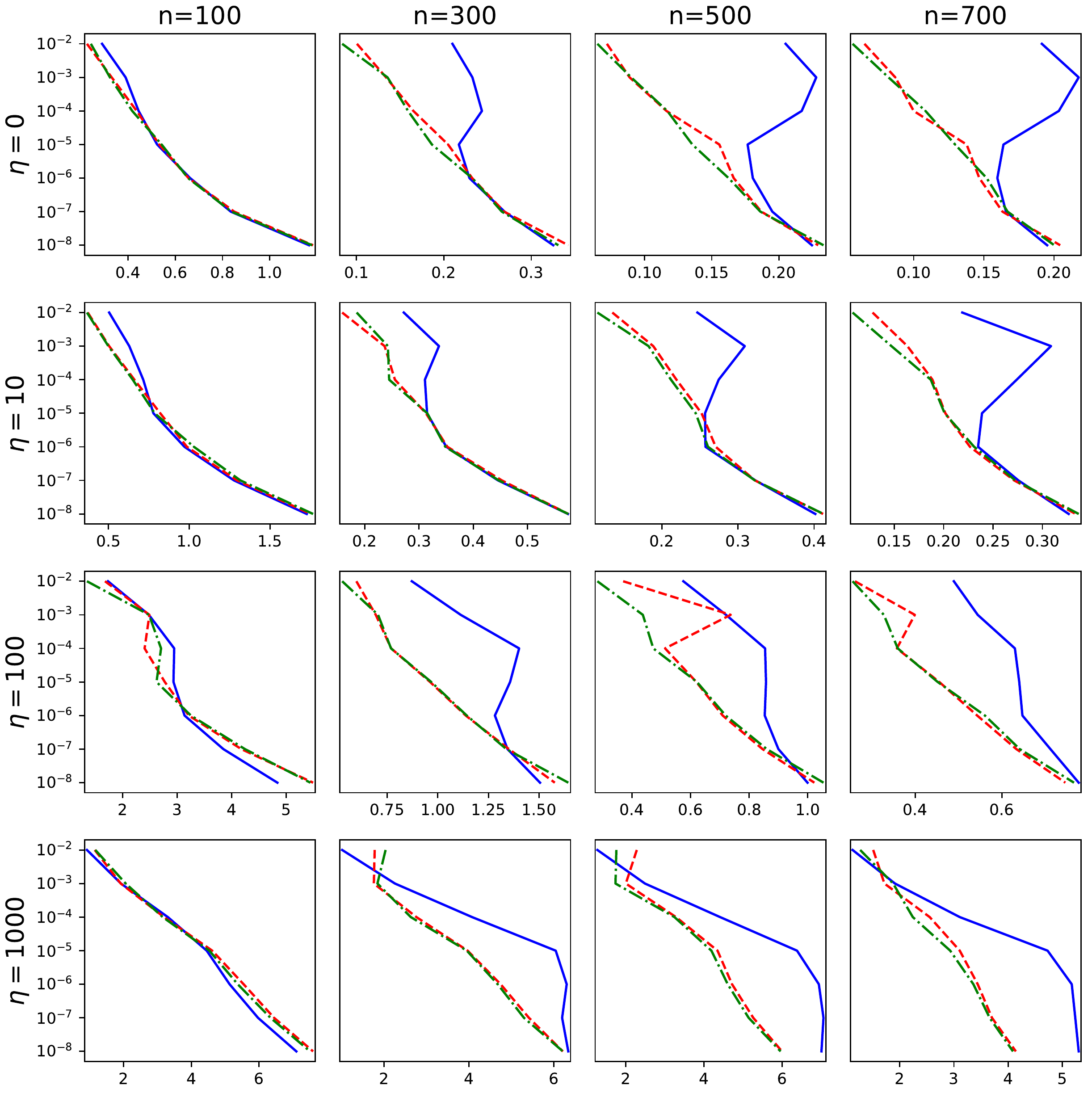}
\par\end{centering}
\caption{The number of normalized Krylov iterations employed by the SDIRK54
scheme (i.e.~computational cost; on the $x$-axis) for a given tolerance
(on the $y$-axis) is shown for equation (\ref{eq:diffadv}) (linear
diffusion-advection equation). The proposed step size controllers
are shown as dashed red lines (penalized variant) and dash-dotted
green lines (non-penalized variant), while the traditional step size
controller is shown in solid blue. The grid size ($n$ is the number
of grid points) and the Péclet number $\eta$ are varied.\label{fig:diffadv-sdirk54}}
\end{figure}

Similar simulations have been conducted for the SDIRK23 and SDIRK54
schemes. These are shown in Figures \ref{fig:diffadv-sdirk23} and
\ref{fig:diffadv-sdirk54}, respectively. The results are similar
to the one discussed in detail for the Crank\textendash Nicolson method.
Although we also see that for (probably unrealistically) high Péclet
numbers the traditional step size controller can slightly outperform
the proposed method. Nevertheless, overall the results reinforce the
significant advantage in performance the proposed method provides.

\section{Nonlinear problems\label{sec:Nonlinear-problems}}

While the linear example in the previous section illustrates very
well the issues with standard step size controllers, ultimately we
are not primarily interested in solving linear diffusion-advection
equations. In addition, an integrator optimized for a rather restricted
class of linear problems is not of much utility in practice. Thus,
in this section we will consider a number of nonlinear problems to
demonstrate that the proposed step size controller also works well
in this regime.

\subsection{Burgers' equation with a reaction term}

As our first example we consider
\begin{equation}
\partial_{t}u(t,x)=\eta u(t,x)\partial_{x}u(t,x)+g(u(t,x)),\label{eq:burger-reaction}
\end{equation}
where we have chosen $g(u)=10(u-2)\sqrt{\vert u-1\vert}$. As usual
periodic boundary conditions on $[0,1]$ are used and the following
initial value is imposed
\[
u(0,x)=2+\epsilon_{1}\sin(\omega_{1}x)+\epsilon_{2}\sin(\omega_{2}x+\varphi)
\]
with $\epsilon_{1}=\epsilon_{2}=10^{-2}$, $\omega_{1}=2\pi$, $\omega_{2}=8\pi$,
and $\varphi=0.3$. The nonlinear reaction satisfies $g(u)<0$ for
$u<2$ and $g(u)>0$ for $u>0$. Thus, the perturbation introduced
in the initial value results in parts of the solution being pulled
towards $1$ while other regions show a growth behavior. The Burgers'
nonlinear, the strength of which is measured by $\eta$, steepens
the gradients between these regions. The integration is performed
until final time $t=0.05$. 

\begin{figure}
\begin{centering}
\includegraphics[width=16cm]{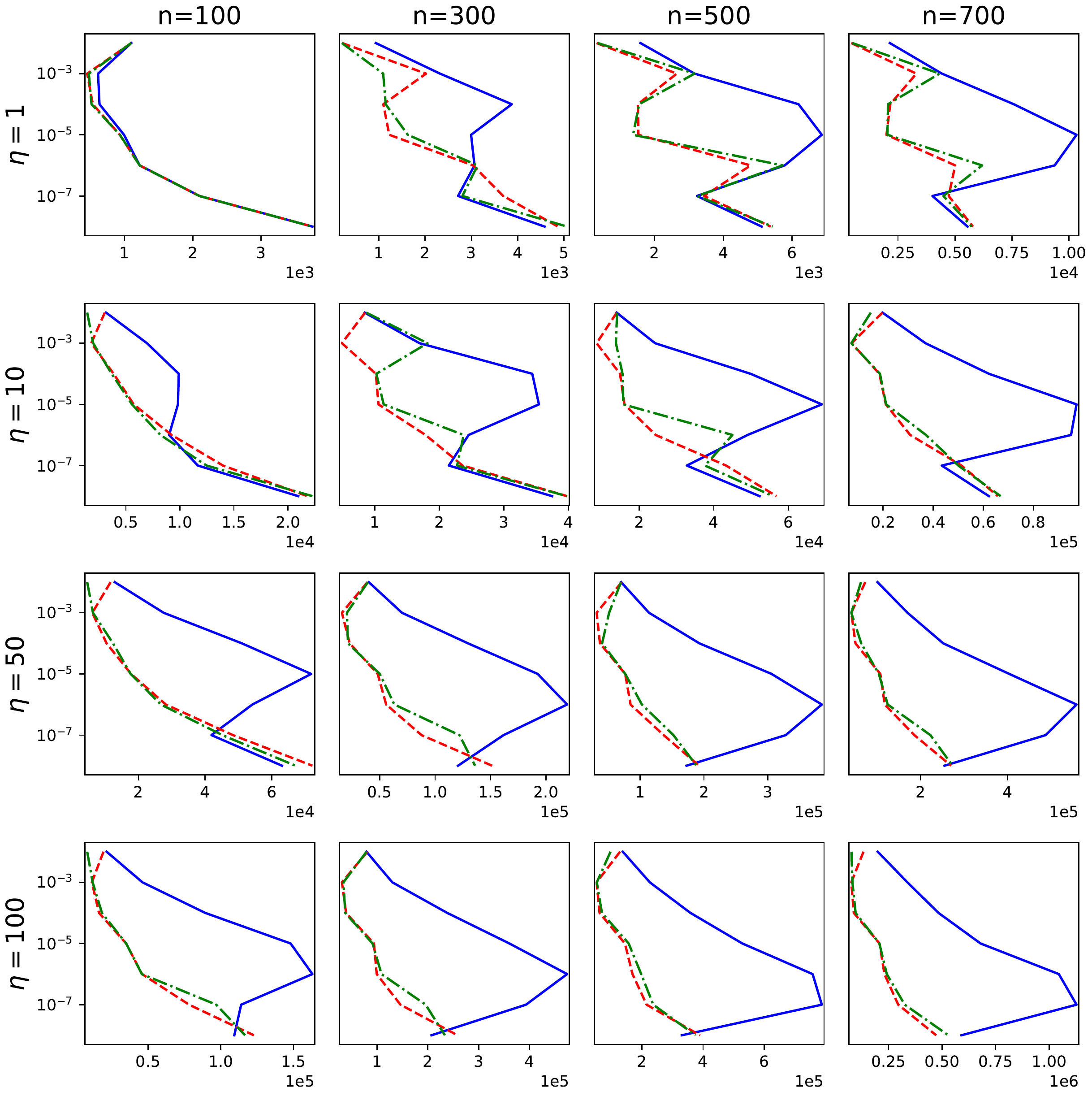}
\par\end{centering}
\caption{The number of Krylov iterations employed by the Crank\textendash Nicolson
scheme (i.e.~computational cost; on the $x$-axis) for a given tolerance
(on the $y$-axis) is shown for equation (\ref{eq:burger-reaction})
(Burgers' equation with a reaction term). The proposed step size controllers
are shown as dashed red lines (penalized variant) and dash-dotted
green lines (non-penalized variant), while the traditional step size
controller is shown in solid blue. The grid size ($n$ is the number
of grid points) and the strength of the Burgers' nonlinearity $\eta$
are varied. \label{fig:instab-cn}}
\end{figure}

\begin{figure}
\begin{centering}
\includegraphics[width=16cm]{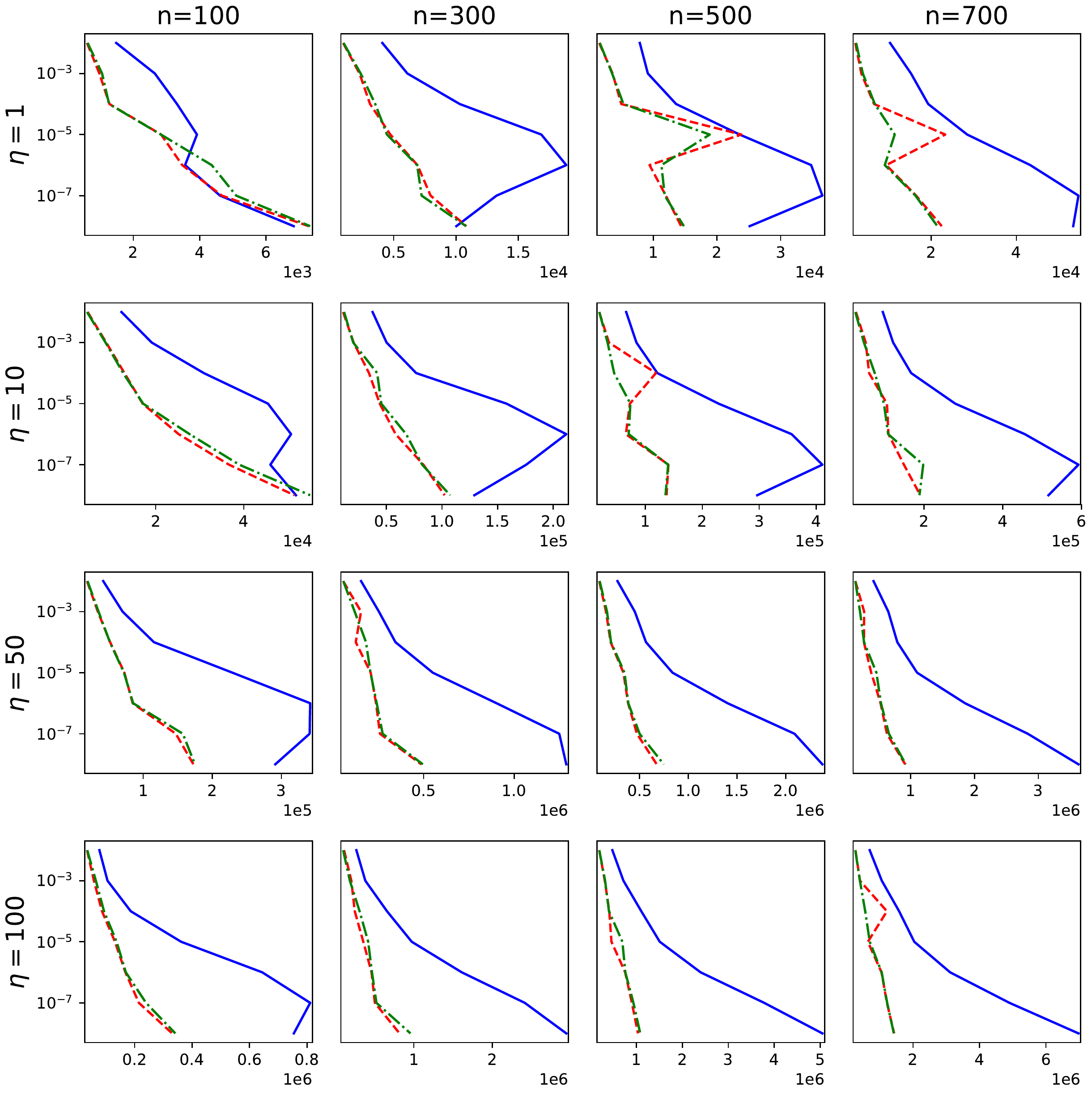}
\par\end{centering}
\caption{The number of Krylov iterations employed by the SDIRK23 scheme (i.e.~computational
cost; on the $x$-axis) for a given tolerance (on the $y$-axis) is
shown for equation (\ref{eq:burger-reaction}) (Burgers' equation
with a reaction term). The proposed step size controllers are shown
as dashed red lines (penalized variant) and dash-dotted green lines
(non-penalized variant), while the traditional step size controller
is shown in solid blue. The grid size ($n$ is the number of grid
points) and the strength of the Burgers' nonlinearity $\eta$ are
varied.\label{fig:instab-sdirk23}}
\end{figure}

\begin{figure}
\begin{centering}
\includegraphics[width=16cm]{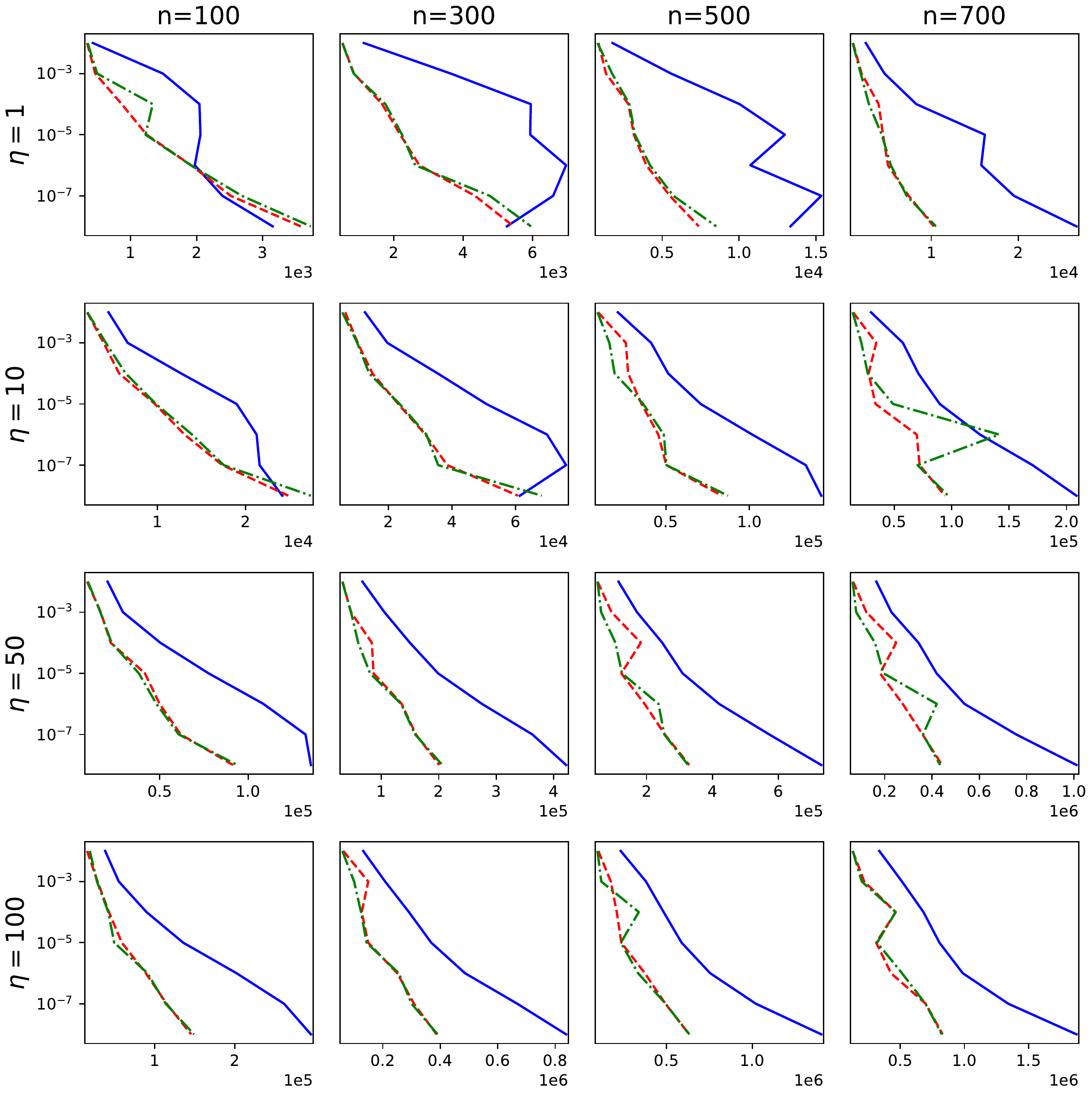}
\par\end{centering}
\caption{The number of Krylov iterations employed by the SDIRK54 scheme (i.e.~computational
cost; on the $x$-axis) for a given tolerance (on the $y$-axis) is
shown for equation (\ref{eq:burger-reaction}) (Burgers' equation
with a reaction term). The proposed step size controllers are shown
as dashed red lines (penalized variant) and dash-dotted green lines
(non-penalized variant), while the traditional step size controller
is shown in solid blue. The grid size ($n$ is the number of grid
points) and the strength of the Burgers' nonlinearity $\eta$ are
varied. \label{fig:instab-sdirk45}}
\end{figure}

The work-precision diagrams for the Crank\textendash Nicolson, SDIRK23,
and SDIRK54 method are shown in Figures \ref{fig:instab-cn}, \ref{fig:instab-sdirk23},
and \ref{fig:instab-sdirk45}, respectively. Overall, we see a significant
improvement compared to the classic step size controller for all numerical
methods and virtually all tolerances. The speedups observed range
up to a factor of five and are most pronounced for medium to large
tolerances.

\begin{figure}
\begin{centering}
\includegraphics[width=16cm]{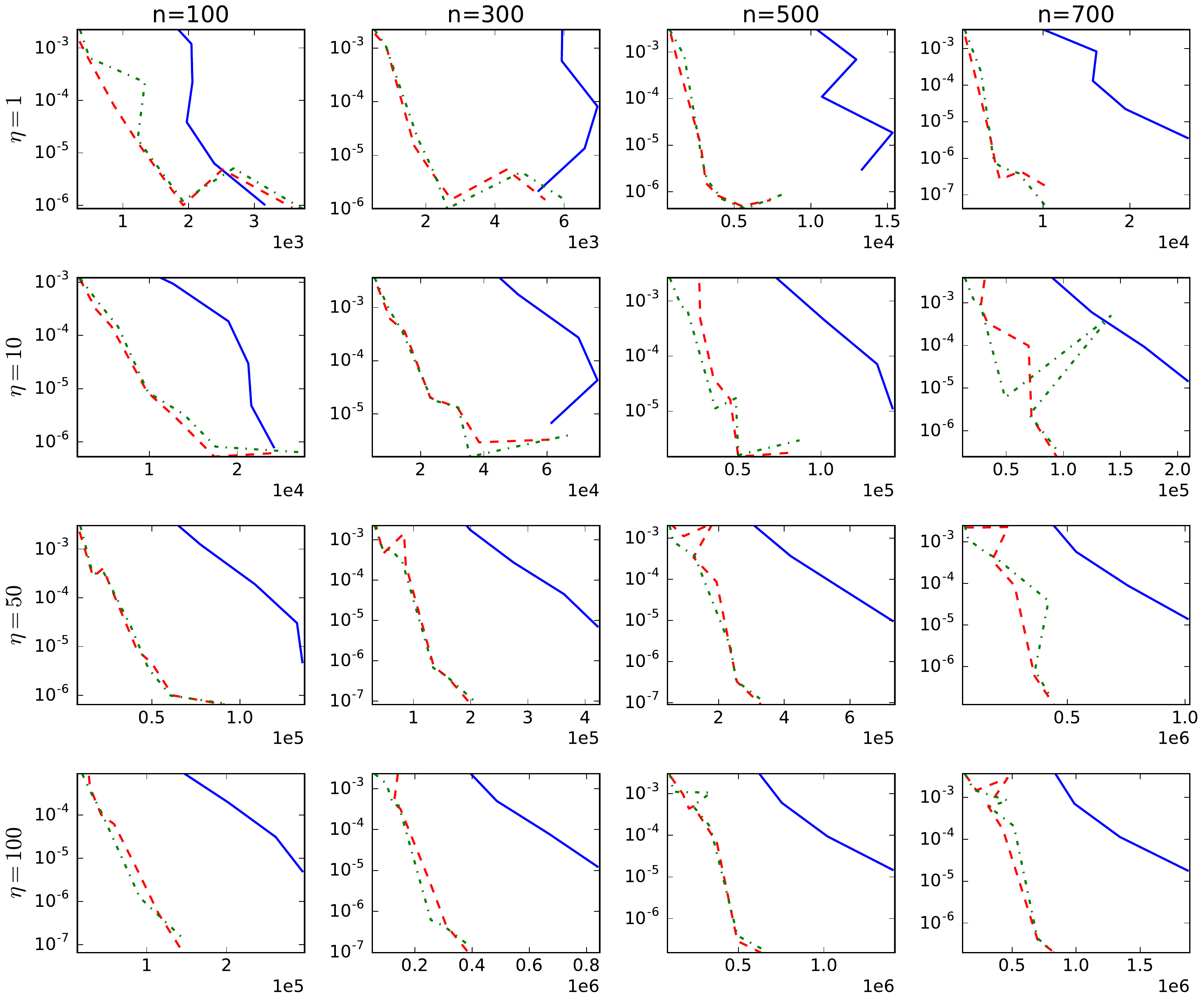}
\par\end{centering}
\caption{The number of Krylov iterations employed by the SDIRK54 scheme (i.e.~computational
cost; on the $x$-axis) for the achieved accuracy at final time $t=0.05$
(on the $y$-axis) is shown for equation (\ref{eq:burger-reaction})
(Burgers' equation with a reaction term). The proposed step size controllers
are shown as dashed red lines (penalized variant) and dash-dotted
green lines (non-penalized variant), while the traditional step size
controller is shown in solid blue. The grid size ($n$ is the number
of grid points) and the strength of the Burgers' nonlinearity $\eta$
are varied. \label{fig:instab-cn-acc}}
\end{figure}

In addition, in Figure \ref{fig:instab-cn-acc} we show the achieved
global error at a final time as a function of the number of Krylov
iterations (this time for the SDIRK54 scheme). As expected, the advantage
of the proposed approach compared to the traditional step size controller
is even more pronounced according to this metric.\textcolor{blue}{{}
}

\subsection{Porous medium equation}

The second example is a porous media equation
\begin{equation}
\partial_{t}u(t,x)=\partial_{xx}(u(t,x)^{m})+\eta\partial_{x}u(t,x),\label{eq:porous}
\end{equation}
where we have chosen $m=2$. Periodic boundary conditions are imposed
on $[0,1]$ and the initial value is given by by a rectangle 
\[
u(0,x)=1+H(x_{1}-x)+H(x-x_{2})
\]
with $x_{1}=0.25$ and $x_{2}=0.6$. Note that in this problem we
have a nonlinear diffusion, we could also write $\partial_{xx}(u^{m})=\partial_{x}(mu^{m-1}\partial_{x}u)$),
coupled to a linear advection. This results in a solution that is
progressively more and more smooth even though we start from a discontinuous
initial value. The equation is integrated until final time $t=10^{-3}$.

\begin{figure}
\begin{centering}
\includegraphics[width=16cm]{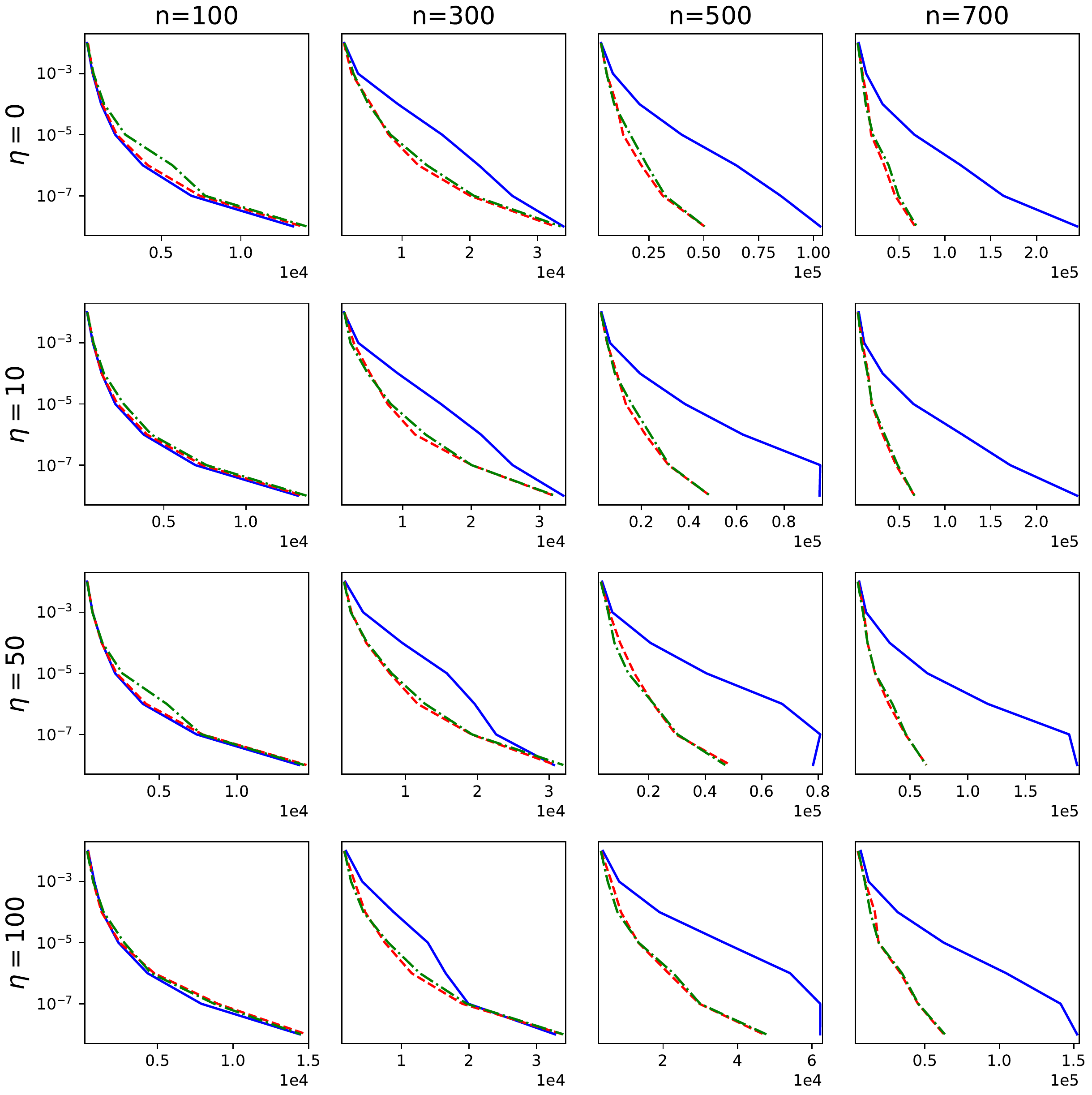}
\par\end{centering}
\caption{The number of Krylov iterations employed by the Crank\textendash Nicolson
scheme (i.e.~computational cost; on the $x$-axis) for a given tolerance
(on the $y$-axis) is shown for equation (\ref{eq:porous}) (Porous
medium equation). The proposed step size controllers are shown as
dashed red lines (penalized variant) and dash-dotted green lines (non-penalized
variant), while the traditional step size controller is shown in solid
blue. The grid size ($n$ is the number of grid points) and the speed
of advection $\eta$ are varied. \label{fig:porous-cn}}
\end{figure}

\begin{figure}
\begin{centering}
\includegraphics[width=16cm]{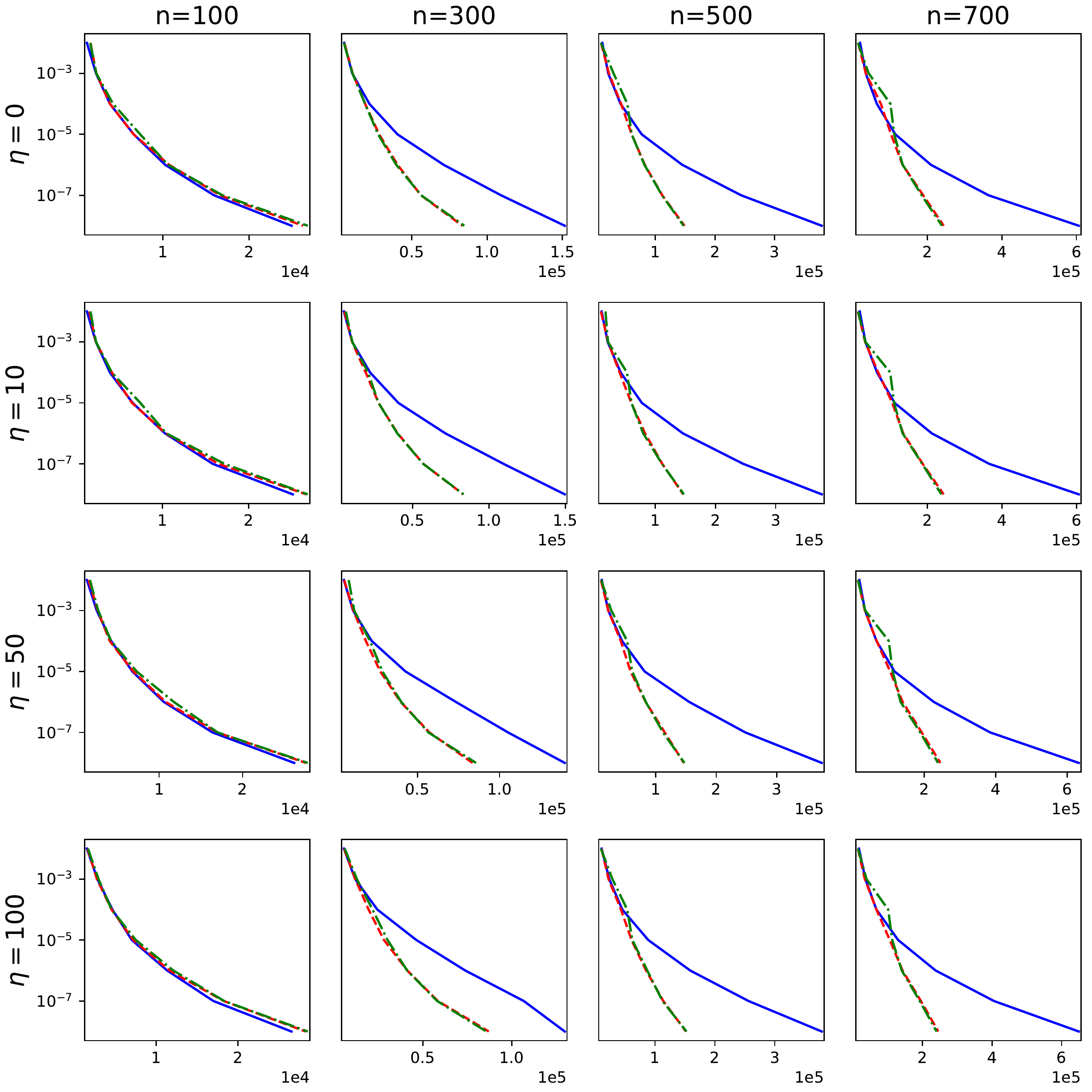}
\par\end{centering}
\caption{The number of Krylov iterations employed by the SDIRK23 scheme (i.e.~computational
cost; on the $x$-axis) for a given tolerance (on the $y$-axis) is
shown for equation (\ref{eq:porous}) (Porous medium equation). The
proposed step size controllers are shown as dashed red lines (penalized
variant) and dash-dotted green lines (non-penalized variant), while
the traditional step size controller is shown in solid blue. The grid
size ($n$ is the number of grid points) and the speed of advection
$\eta$ are varied.\label{fig:porous-sdirk23}}
\end{figure}

\begin{figure}
\begin{centering}
\includegraphics[width=16cm]{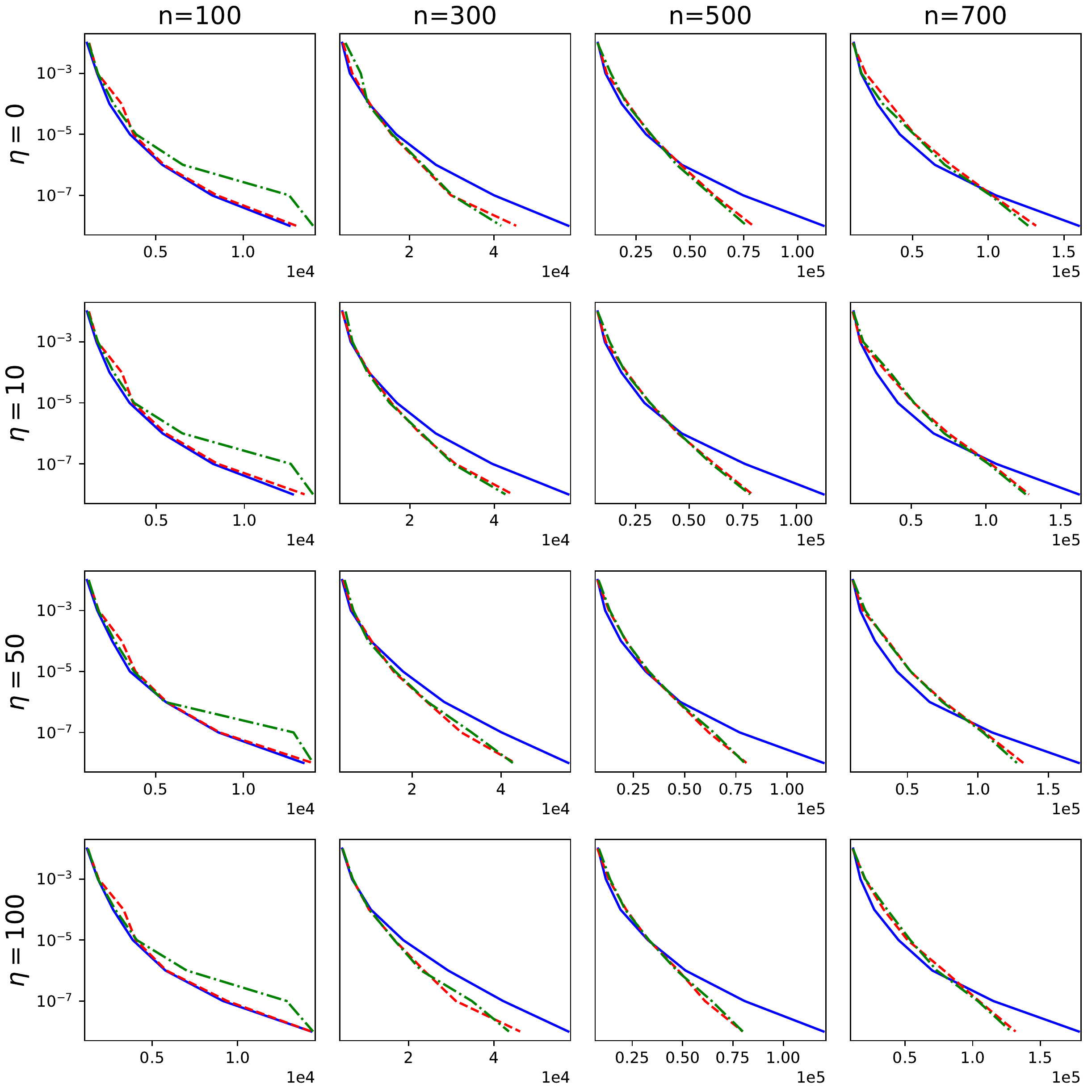}
\par\end{centering}
\caption{The number of Krylov iterations employed by the SDIRK54 scheme (i.e.~computational
cost; on the $x$-axis) for a given tolerance (on the $y$-axis) is
shown for equation (\ref{eq:porous}) (Porous medium equation). The
proposed step size controllers are shown as dashed red lines (penalized
variant) and dash-dotted green lines (non-penalized variant), while
the traditional step size controller is shown in solid blue. The grid
size ($n$ is the number of grid points) and the speed of advection
$\eta$ are varied.\label{fig:porous-sdirk45}}
\end{figure}

The work-precision diagrams for the Crank\textendash Nicolson, SDIRK23,
and SDIRK54 method are shown in Figures \ref{fig:porous-cn}, \ref{fig:porous-sdirk23},
and \ref{fig:porous-sdirk45}, respectively. For the Crank\textendash Nicolson
method and SDIRK23 we see a significant increase in performance, in
particular, as we use more grid points to discretize the problem.
The maximal achieved speedups for the numerical simulation shown is
approximately a factor of four. The gains are significantly more modest
for the SDIRK54 method. In this case maximal gains are on the order
of 50\% and they only manifest themselves at tolerances well below
$10^{-6}$. For low tolerances the proposed step size controller can
be slightly slower than the traditional approach. It is in this regime
that the penalized variant performs significantly better compared
to the non-penalized variant (otherwise the two variants show similar
performance). Nevertheless, overall a clear improvement in performance
can still be observed.

\subsection{Viscous Burgers' equation}

The third example is the viscous Burgers' equation
\begin{equation}
\partial_{t}u(t,x)=\partial_{xx}u(t,x)-\eta u(t,x)\partial_{x}u(t,x).\label{eq:viscous-burger}
\end{equation}
Periodic boundary conditions are imposed on $[0,1]$ and the initial
value
\[
u(0,x)=1+e\exp\left(\frac{-1}{1-(2x-1)^{2}}\right)+\frac{1}{2}\exp\left(-\frac{(x-x_{0})^{2}}{2\sigma^{2}}\right),
\]
where $x_{0}=0.9$ and $\sigma=0.02$ is chosen. This corresponds
to a $\mathcal{C}^{\infty}$ bump at the tail of which a Gaussian
of smaller amplitude is added. In this example the nonlinearity tries
to create strong gradients, while the diffusion counteracts that.
That is, $\eta$ is a measure of how strong a resistance is provided
to the homogenization of the solution. The problem is integrated to
final time $t=10^{-2}$.

\begin{figure}
\begin{centering}
\includegraphics[width=16cm]{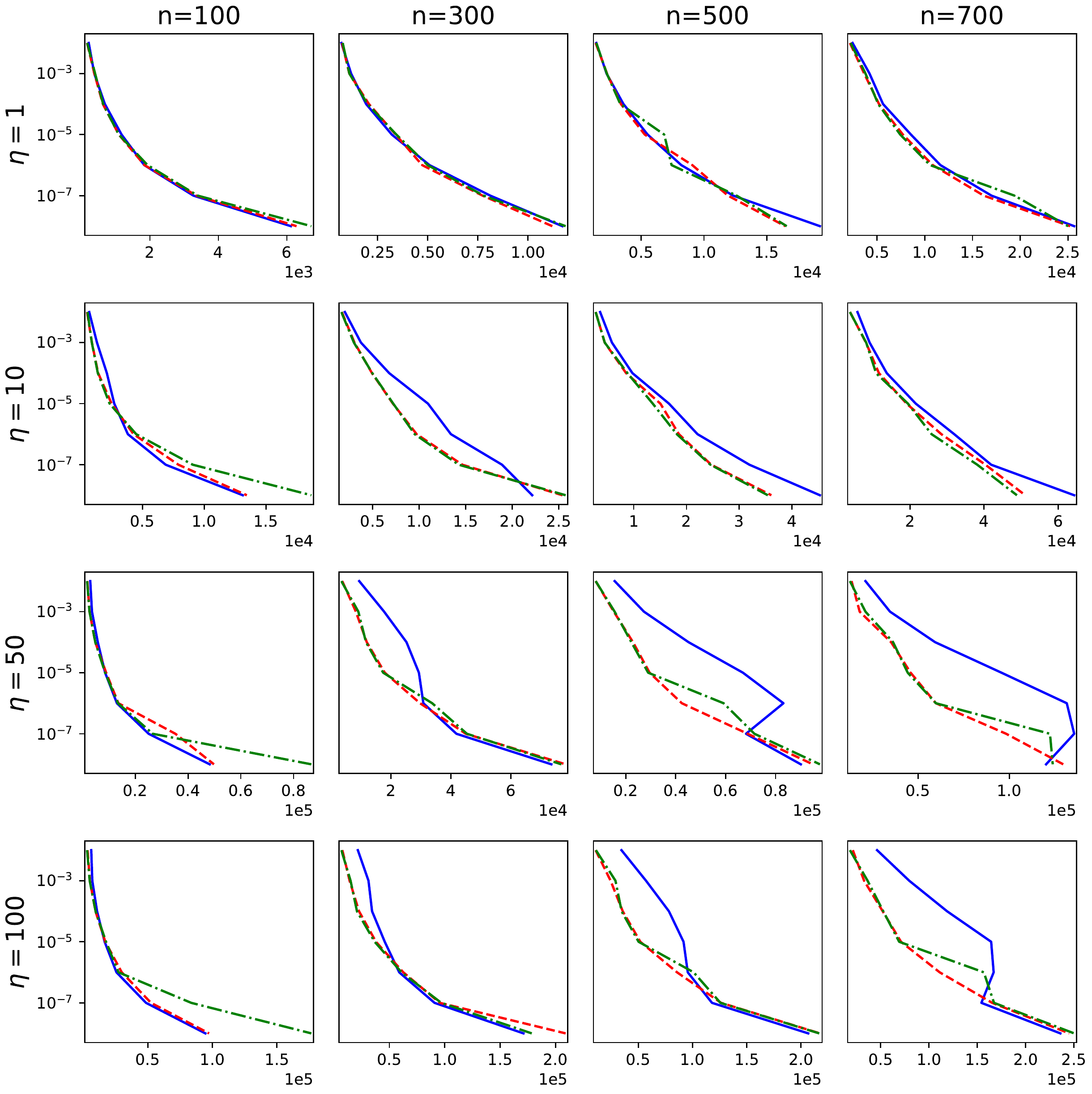}
\par\end{centering}
\caption{The number of Krylov iterations employed by the Crank\textendash Nicolson
scheme (i.e.~computational cost; on the $x$-axis) for a given tolerance
(on the $y$-axis) is shown for equation (\ref{eq:viscous-burger})
(viscous Burgers' equation). The proposed step size controllers are
shown as dashed red lines (penalized variant) and dash-dotted green
lines (non-penalized variant), while the traditional step size controller
is shown in solid blue. The grid size ($n$ is the number of grid
points) and the strength of the nonlinear advection $\eta$ are varied.\label{fig:burger-cn}}
\end{figure}

\begin{figure}
\begin{centering}
\includegraphics[width=16cm]{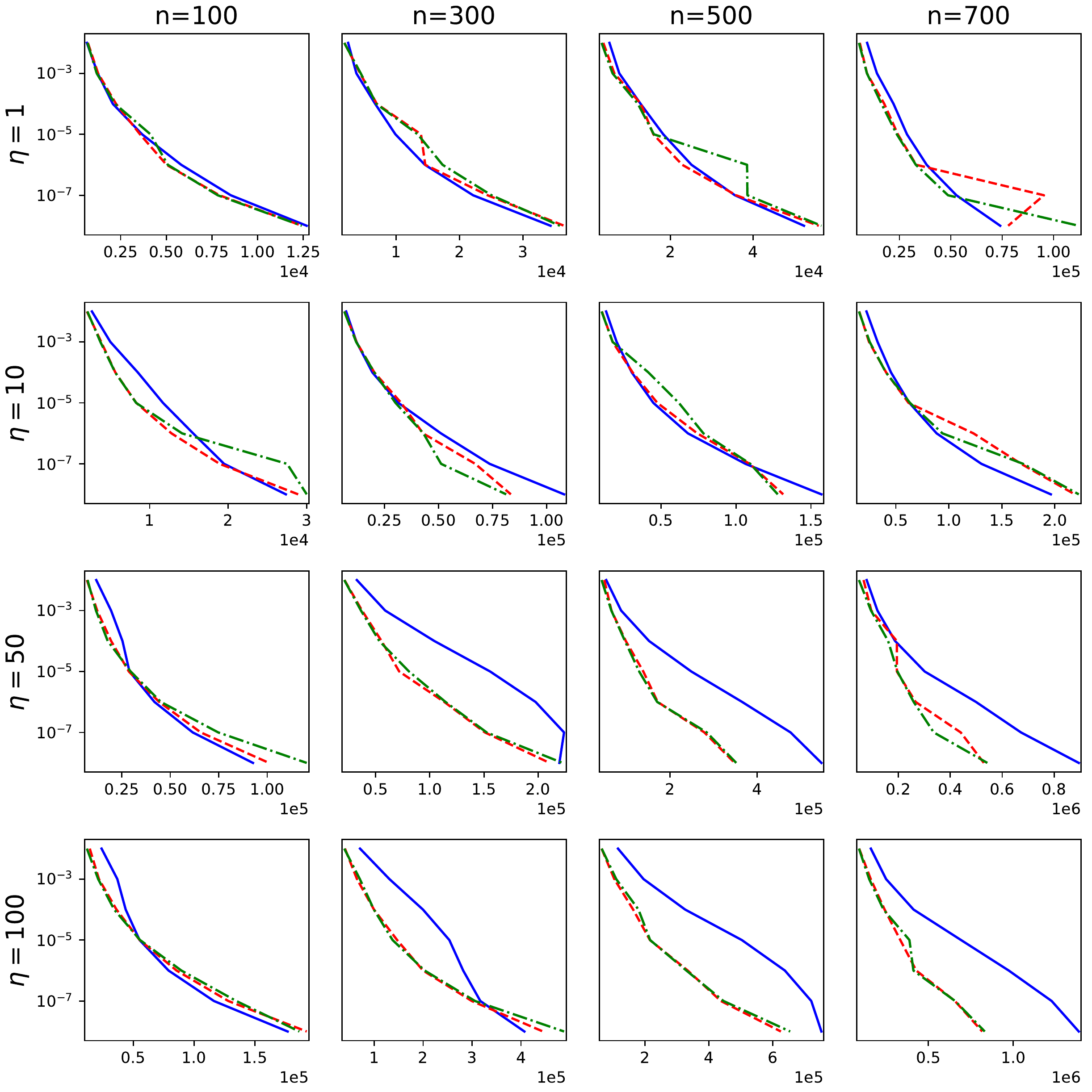}
\par\end{centering}
\caption{The number of Krylov iterations employed by the SDIRK23 scheme (i.e.~computational
cost; on the $x$-axis) for a given tolerance (on the $y$-axis) is
shown for equation (\ref{eq:viscous-burger}) (viscous Burgers' equation).
The proposed step size controllers are shown as dashed red lines (penalized
variant) and dash-dotted green lines (non-penalized variant), while
the traditional step size controller is shown in solid blue. The grid
size ($n$ is the number of grid points) and the strength of the nonlinear
advection $\eta$ are varied.\label{fig:burger-sdirk23}}
\end{figure}

\begin{figure}
\begin{centering}
\includegraphics[width=16cm]{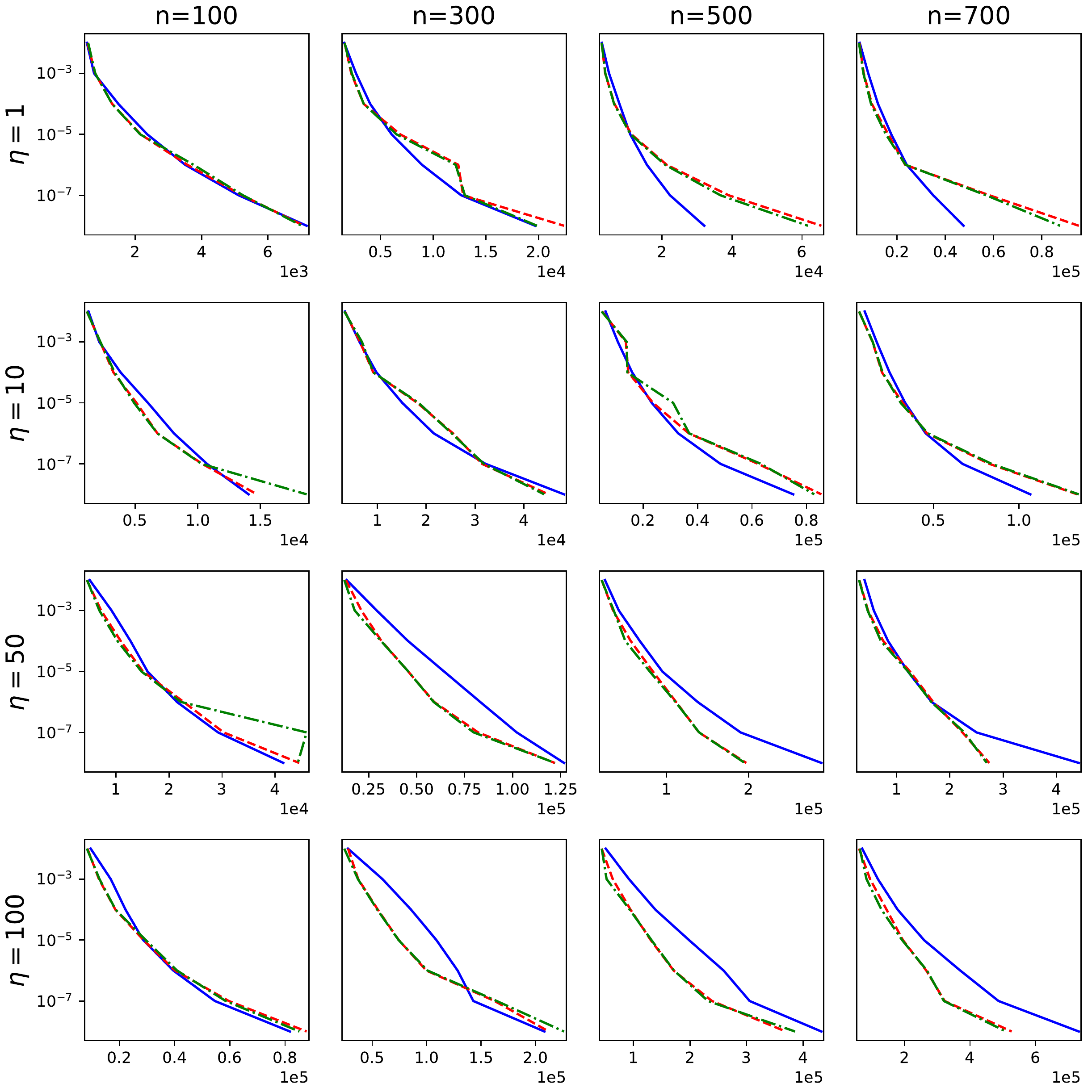}
\par\end{centering}
\caption{The number of Krylov iterations employed by the SDIRK54 scheme (i.e.~computational
cost; on the $x$-axis) for a given tolerance (on the $y$-axis) is
shown for equation (\ref{eq:viscous-burger}) (viscous Burgers' equation).
The proposed step size controllers are shown as dashed red lines (penalized
variant) and dash-dotted green lines (non-penalized variant), while
the traditional step size controller is shown in solid blue. The grid
size ($n$ is the number of grid points) and the strength of the nonlinear
advection $\eta$ are varied.\label{fig:burger-sdirk45}}
\end{figure}

The work-precision diagrams for the Crank\textendash Nicolson, SDIRK23,
and SDIRK54 method are shown in Figures \ref{fig:burger-cn}, \ref{fig:burger-sdirk23},
and \ref{fig:burger-sdirk45}, respectively. For relatively low $\eta$
only small improvements or even a slowdown can be observed. However,
as we increase $\eta$ the proposed step size controller shows a significant
advantage for all numerical methods. Speedups up to a factor of two
are observed, particular as the number of grid points is increased.

\subsection{Allen\textendash Cahn equation}

As the fourth example we consider the one-dimensional Allen-Cahn equation
\begin{equation}
\partial_{t}u(t,x)=\partial_{xx}u(t,x)+\eta u(t,x)(1-u(t,x)^{2}).\label{eq:allencahn}
\end{equation}
Periodic boundary conditions are imposed on $[0,1]$ and the following
initial value
\[
u(0,x)=A\left(1+\cos\omega_{1}x\right),
\]
with $A=\tfrac{1}{10}$ and $\omega_{1}=2\pi$ is used. This problem
does not include an advection term. As a consequence, the linear part
of the right-hand side is a symmetric matrix. The interesting behavior
of the Allen-Cahn equation is due to the fact that the nonlinear reaction
term pulls the solution to either $0$, 1, or $-1$, while the diffusion
tries to homogenize the solution. This results, for larger $\eta$,
in regions of space that are separated by relatively sharp gradients.
The problem is integrated to final time $t=2\cdot10^{-2}$.

\begin{figure}
\begin{centering}
\includegraphics[width=16cm]{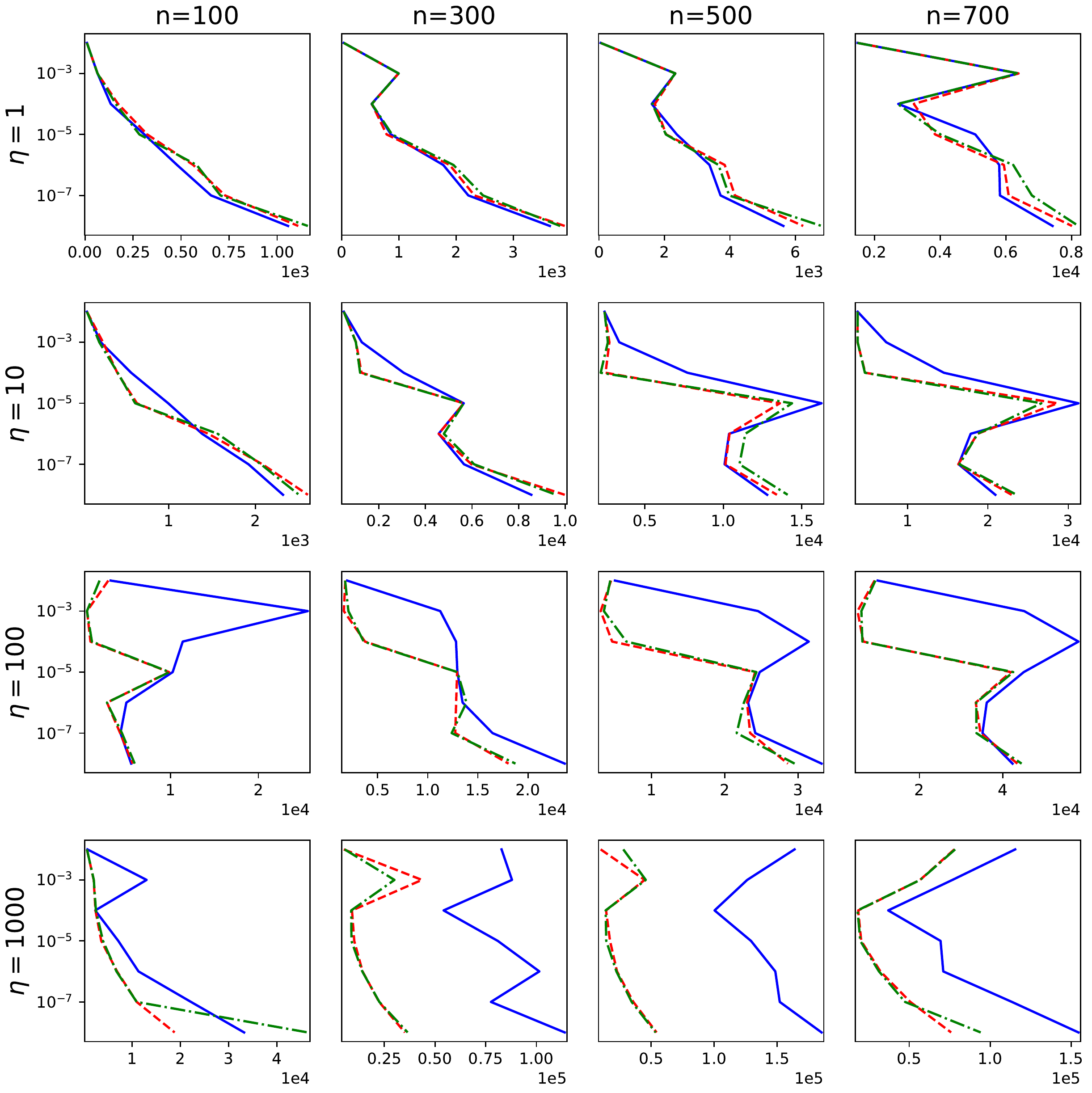}
\par\end{centering}
\caption{The number of Krylov iterations employed by the Crank\textendash Nicolson
scheme (i.e.~computational cost; on the $x$-axis) for a given tolerance
(on the $y$-axis) is shown for equation (\ref{eq:allencahn}) (Allen-Cahn
equation). The proposed step size controllers are shown as dashed
red lines (penalized variant) and dash-dotted green lines (non-penalized
variant), while the traditional step size controller is shown in solid
blue. The grid size ($n$ is the number of grid points) and the strength
of the nonlinear reaction $\eta$ are varied.\label{fig:allencahn-cn}}
\end{figure}

\begin{figure}
\begin{centering}
\includegraphics[width=16cm]{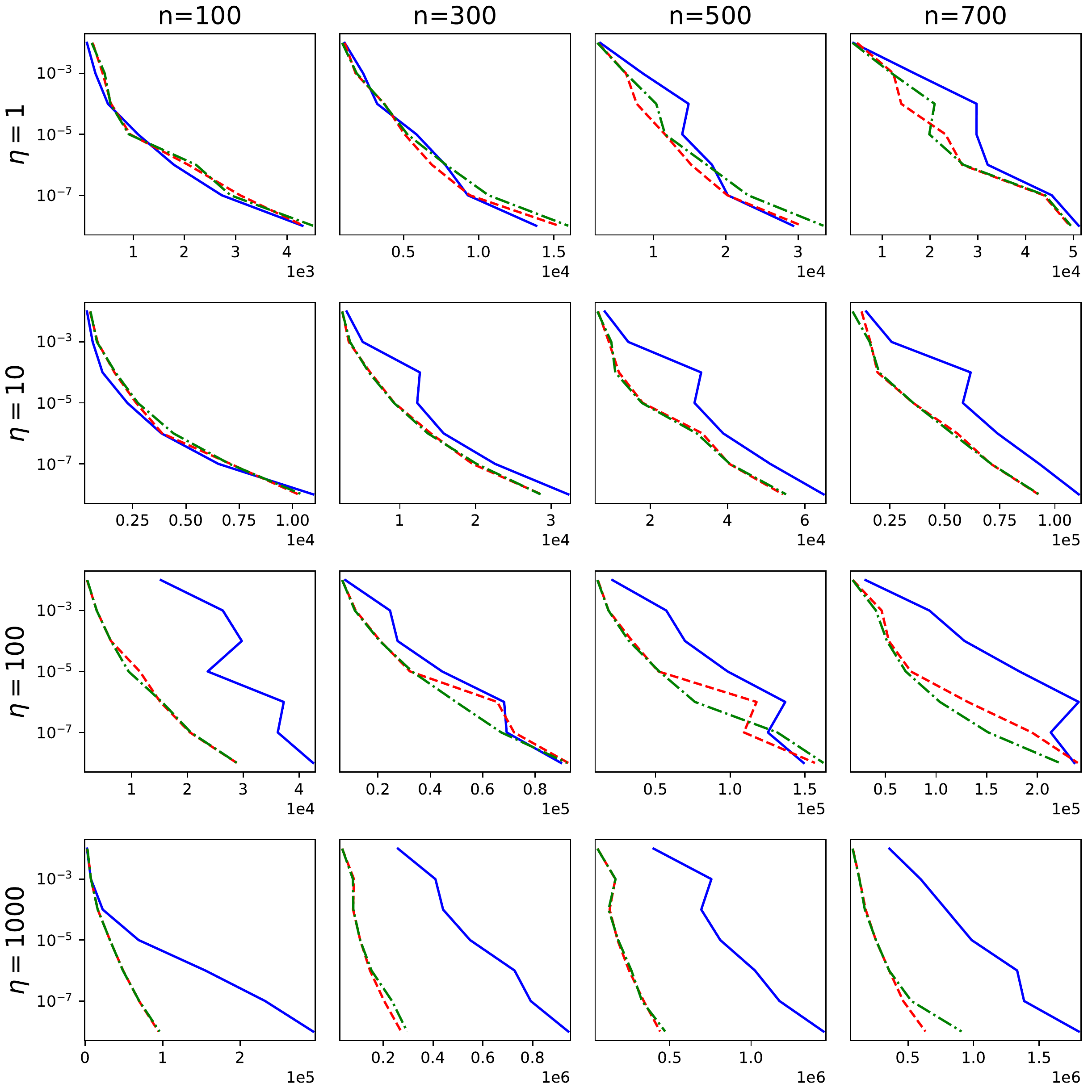}
\par\end{centering}
\caption{The number of Krylov iterations employed by the SDIRK23 scheme (i.e.~computational
cost; on the $x$-axis) for a given tolerance (on the $y$-axis) is
shown for equation (\ref{eq:allencahn}) (Allen-Cahn equation). The
proposed step size controllers are shown as dashed red lines (penalized
variant) and dash-dotted green lines (non-penalized variant), while
the traditional step size controller is shown in solid blue. The grid
size ($n$ is the number of grid points) and the strength of the nonlinear
reaction $\eta$ are varied.\label{fig:allencahn-sdirk23}}
\end{figure}

\begin{figure}
\begin{centering}
\includegraphics[width=16cm]{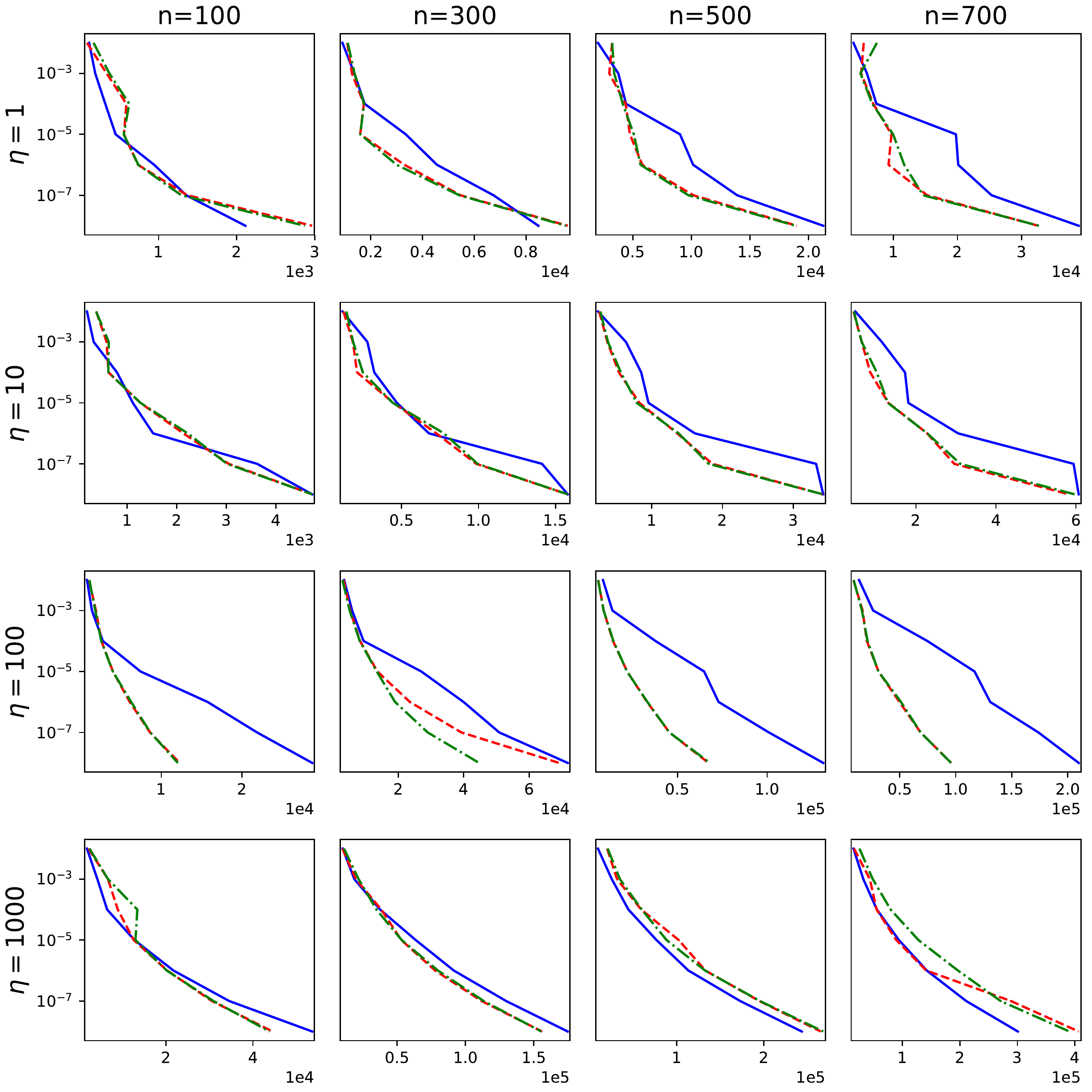}
\par\end{centering}
\caption{The number of Krylov iterations employed by the SDIRK54 scheme (i.e.~computational
cost; on the $x$-axis) for a given tolerance (on the $y$-axis) is
shown for equation (\ref{eq:allencahn}) (Allen-Cahn equation). The
proposed step size controllers are shown as dashed red lines (penalized
variant) and dash-dotted green lines (non-penalized variant), while
the traditional step size controller is shown in solid blue. The grid
size ($n$ is the number of grid points) and the strength of the nonlinear
reaction $\eta$ are varied.\label{fig:allencahn-sdirk45}}
\end{figure}

The work-precision diagrams for the Crank\textendash Nicolson, SDIRK23,
and SDIRK54 method are shown in Figures \ref{fig:allencahn-cn}, \ref{fig:allencahn-sdirk23},
and \ref{fig:allencahn-sdirk45}, respectively. For the Crank\textendash Nicolson
and the SDIRK23 scheme large increases in performance can be observed
for large $\eta$. For the SDIRK54 method the proposed step size controller
shows significant savings in computational effort for $\eta\leq100$
and similar performance for $\eta=1000$. Although there are some
small differences between the penalized and non-penalized variants,
in the present test these two methods perform very similar. 

\section{Brusselator in two dimensions}

As the final example we consider the two-dimensional Brusselator given
by
\begin{align}
\partial_{t}u(t,x,y) & =\alpha\Delta u(t,x,y)+1+u^{2}v-4.4u+f(t,x,y)\nonumber \\
\partial_{t}v(t,x,y) & =\alpha\Delta u(t,x,y)+3.4u-u^{2}v,\label{eq:brusselator}
\end{align}
where $\Delta=\partial_{xx}+\partial_{yy}$ is the Laplacian and $\alpha=0.1$.
Periodic boundary conditions on $[0,1]^{2}$ are imposed and the following
initial value
\[
u(0,x,y)=22y(1-y)^{3/2},\qquad v(0,x,y)=27x(1-x)^{3/2}
\]
is selected. The function $f$ is a source term and is chosen such
that $f(t,x,y)=5$ if $(x-0.3)^{2}+(y-0.6)^{2}\leq0.1^{2}$ and $t\geq1.1$.
Otherwise $f(t,x,y)$ is set to zero. The problem is integrated to
final time $t=11.5$. This is the problem considered in \cite[p. 151--152]{hairerII}.

\begin{figure}
\begin{centering}
Crank\textendash Nicolson\smallskip{}
\par\end{centering}
\begin{centering}
\includegraphics[width=16cm]{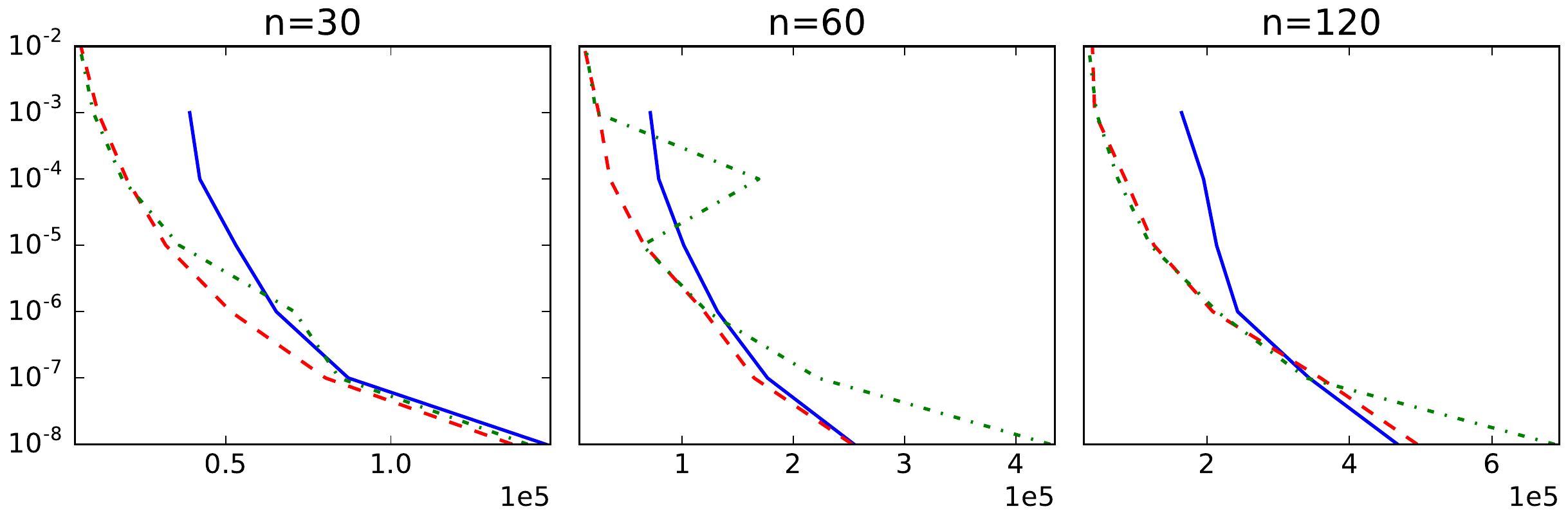}
\par\end{centering}
\begin{centering}
SDIRK23
\par\end{centering}
\begin{centering}
\smallskip{}
\par\end{centering}
\begin{centering}
\includegraphics[width=16cm]{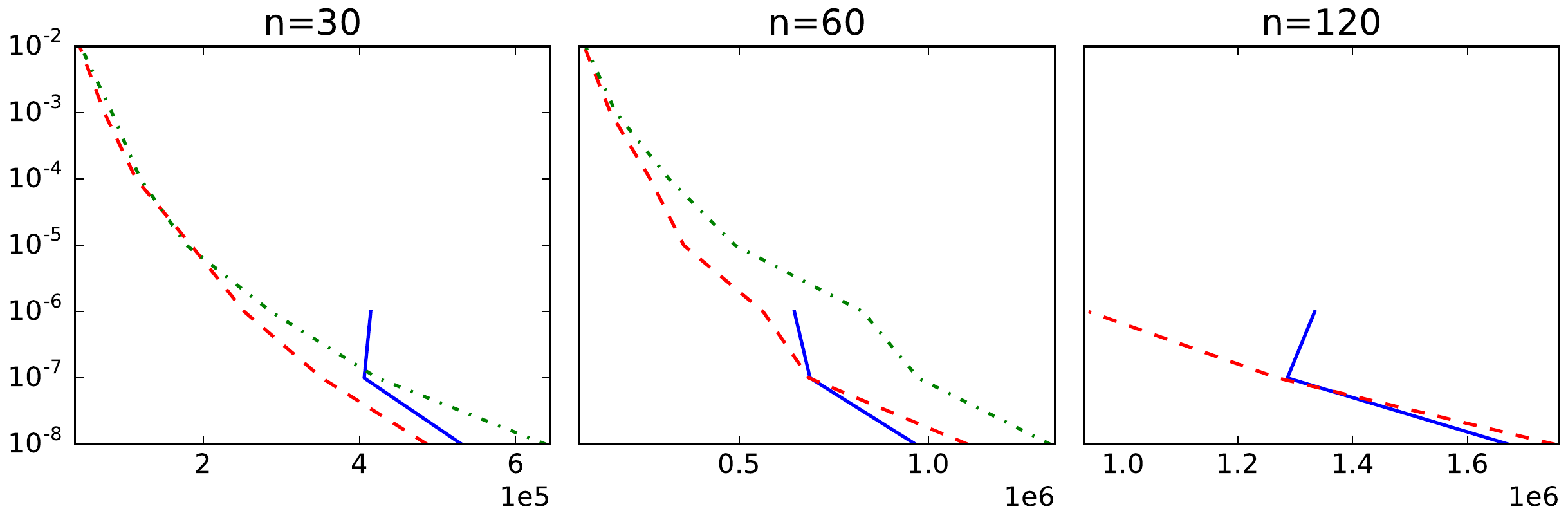}
\par\end{centering}
\begin{centering}
SDIRK54\smallskip{}
\par\end{centering}
\begin{centering}
\includegraphics[width=16cm]{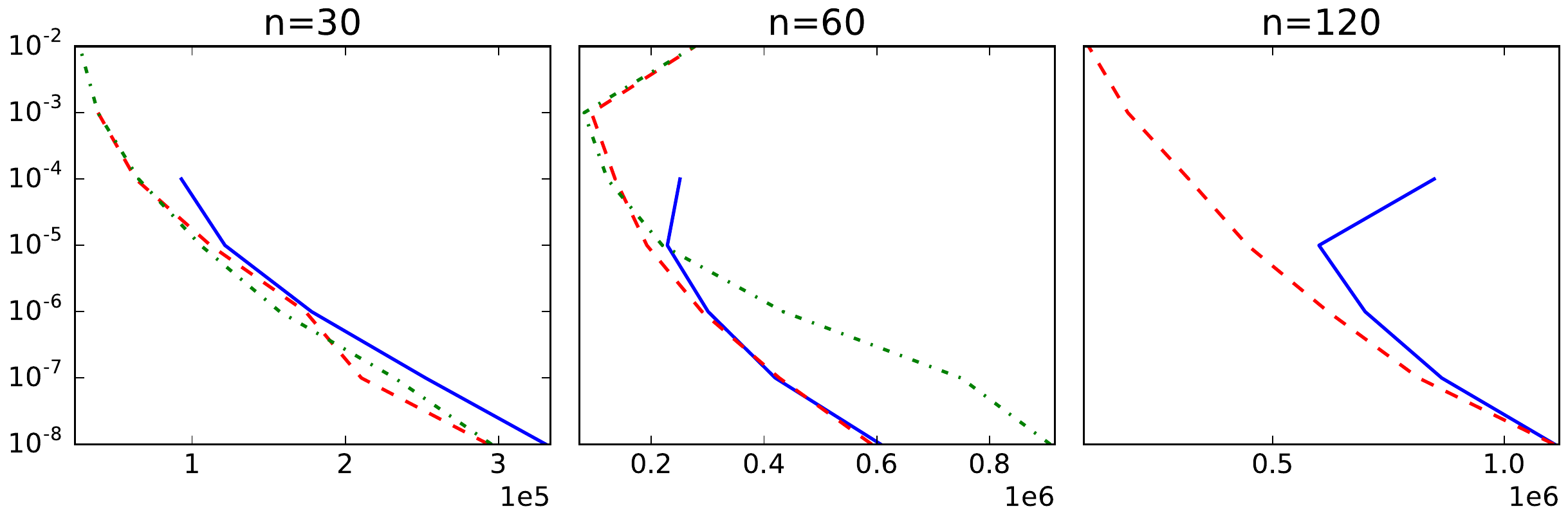}
\par\end{centering}
\caption{The number of Krylov iterations employed by the Crank\textendash Nicolson/SDIRK23/SDIRK54
schemes (i.e.~computational cost; on the $x$-axis) for a given tolerance
(on the $y$-axis) is shown for equation (\ref{eq:brusselator}) (the
two-dimensional Brusselator). The proposed step size controllers are
shown as dashed red lines (penalized variant) and dash-dotted green
lines (non-penalized variant), while the traditional step size controller
is shown in solid blue. The grid size ($n$ is the number of grid
points per direction) is varied.\label{fig:brusselator}}
\end{figure}

The work-precision diagrams for the Crank\textendash Nicolson, SDIRK23,
and SDIRK54 methods are shown in Figure \ref{fig:brusselator}. The
penalized variant of the proposed step size controller shows superior
performance in virtually all configuration. The non-penalized variants
performs worse. This is particularly true for the finest grid with
the SDIRK23 and SDIRK54 schemes, where it does not converge in a reasonable
amount of time. However, in the other configurations it is still able
to significantly outperform the traditional step size controller in
the low to medium precision regime. The penalized version of the proposed
step size controller mostly avoids the inverse C curve and is more
robust. The latter is most apparent for the SDIRK23 scheme, where
the standard step size controller, for tolerances above $10^{-6}$,
does not produce a solution within a reasonable amount of time (note
that this issue was also reported in the context of CVODE applied
to a magnetohydrodynamics problem \cite{einkemmer2017}).

\section{Conclusion\label{sec:Conclusion}}

We have demonstrated that the proposed adaptive step size selection
strategy result in significant improvements compared to more traditional
approaches in the context of a number of one-dimensional test problems.
Speedups of up to a factor of five have been observed and significant
increases in performance are seen in almost all problems. In addition,
the inverse C curve is straightened out in almost all configurations
which makes the step size controller more predictable in practice.
In addition, we have considered the two-dimensional Brusselator, where
similar conclusions can be drawn for the penalized variant. 

It is also interesting to note that the speedup observed is most pronounced
for the Crank\textendash Nicolson method, which is still widely used
by physicists and engineers. This might be considered a further disadvantage
of this method (i.e.~traditional step size controllers do not work
as well). Whether this is generally true for methods that fail to
be L-stable might warrant further investigation.

The method proposed here is relatively simple since it only requires
a limited set of parameters. Consequently, optimizing for these parameters
is relatively straightforward and our approach, which only uses selected
samples of a linear problem, demonstrates that such an approach can
be successfully generalized to different linear and even nonlinear
problems. However, one could envisage an approach that makes use of
data obtained in previous time steps as well as allows for more general
mappings. To train such a model, however, would require an extensive
set of representative test problems. We consider this as future work.
In addition, our plan is to extend the step size controller proposed
here to more elaborate implicit methods (such as Radau and Gauss methods)
and consider physically more realistic problems in multiple dimensions.

\bibliographystyle{plain}
\bibliography{paper-stepsizectrl}

\end{document}